\documentclass[lettersize,journal]{IEEEtran}
\usepackage{amsmath,amsfonts}
\usepackage{algorithmic}
\usepackage{algorithm}
\usepackage{array}
\usepackage[caption=false,font=normalsize,labelfont=sf,textfont=sf]{subfig}
\usepackage{textcomp}
\usepackage{stfloats}
\usepackage{url}
\usepackage{verbatim}
\usepackage{graphicx}
\usepackage{cite}

\usepackage{amsmath,amssymb,amsfonts}
\usepackage{graphicx}
\usepackage{textcomp}
\usepackage{mystyle}

\begin{document}

\title{Algorithms based on operator-averaged operators}

\author{Miguel~Sim\~{o}es
\thanks{M. Sim\~{o}es is with the Dept.\ Electr.\ Eng.\ (ESAT) -- STADIUS, KU Leuven, Leuven, Belgium; e-mail: miguel.simoes@kuleuven.be.}
}

\markboth{Preprint}%
{A Sample Article Using IEEEtran.cls for IEEE Journals}

\maketitle

\begin{abstract}
A class of algorithms comprised by certain semismooth Newton and active-set methods is able to solve convex minimization problems involving sparsity-inducing regularizers very rapidly; the speed advantage of methods from this class is a consequence of their ability to benefit from the sparsity of the corresponding solutions by solving smaller inner problems than conventional methods. 
The convergence properties of such conventional methods (e.g., the forward--backward and the proximal--Newton ones) can be studied very elegantly under the framework of iterations of scalar-averaged operators---this is not the case for the aforementioned class. 
However, we show in this work that by instead considering operator-averaged operators, one can indeed study methods of that class, and also to derive algorithms outside of it that may be more convenient to implement than existing ones. 
Additionally, we present experiments whose results suggest that methods based on operator-averaged operators achieve substantially faster convergence than conventional ones.
\end{abstract}


\section{Introduction}
\IEEEPARstart{M}{any} large-scale inverse problems in signal and image processing can be formulated as minimization problems whose objective functions are sums of proper lower-semicontinuous convex functions. 
An example is
\begin{equation} \label{eq:split}
	\underset{\mathbf{x} \in \mathbb{R}^n}{\text{minimize}} \quad f(\mathbf{x}) + g(\mathbf{x}) + \sum_{j=1}^N h_j(\mathbf{L}_j \mathbf{x}),
\end{equation}
where $f: \mathbb{R}^n \to \; ]-\infty,+\infty]$ is assumed to be \emph{smooth}---i.e., to be differentiable with a Lipschitz-continuous gradient---, but $g: \mathbb{R}^n \to \; ]-\infty,+\infty]$ and the $N$ functions $h_j: \mathbb{R}^{m_j} \to \; ]-\infty,+\infty]$, with $\mathbf{L}_j \in \mathbb{R}^{m_j \times n}$ for $j = \{1, \dots, N\}$, may not be. 
Typical examples of problems that fit this formulation are problems where $f$ is a data-fitting term, and $g$ and $h_j$ are regularizers. 
This is the case of the regression method known as \ac{LASSO}, where $f = \norm{\mathbf{H} \cdot - \mathbf{b}}{2}^2$, $g = \mu \norm{\cdot}{1}$, and $N=0$, for given $\mathbf{H} \in \mathbb{R}^{m \times n}$, $\mathbf{b} \in \mathbb{R}^{m}$, and $\mu > 0$. 
The functions $g$ and $h_j$ can also be indicator functions of convex sets, allowing one to consider constrained problems within this framework. 
For example, if in the problem just discussed we additionally consider $h_1$ (with $N=1$) to be the indicator function of the set of nonnegative real numbers, i.e., if $h_1 = \delta_{\mathbb{R}^+_0}$, we are constraining $\mathbf{x} \geq 0$.

Methods known as \emph{splitting methods} convert~\eqref{eq:split} into a sequence of separable subproblems that are easier to solve than the original one. 
One way of studying the convergence properties of such methods is to consider that their iterations can be described by sequences of \emph{averaged operators} (see below for a definition). 
This study is advantageous because it unifies and simplifies the analysis of a large number of algorithms.
However, not all splitting methods can be described as sequences of averaged operators. 
This is the case, for example, of methods belonging to the class of algorithms comprised by semismooth Newton and active-set methods.
Such methods have been applied to develop fast solvers for problems involving sparsity-inducing regularizers, since their iterations involve computing the solution of subproblems whose dimensionality is smaller than in other methods. 
In this work, we consider an extension of the concept of averaged operators, and use this extension to study some of the algorithms belonging to the class just discussed. 
Additionally, we develop new methods that do not belong to this class but that are instances of the aforementioned extension, and that have convenient practical applications.
In the remainder of this section, we list the notation used throughout this work (Subsection~\ref{sec:notation}), provide a brief overview of the \ac{FB} method and some of its variants as to motivate our contributions (Subsection~\ref{sec:background}), and conclude by enumerating those contributions and providing an outline of the paper (Subsection~\ref{sec:outline}).

\subsection{Notation} \label{sec:notation}

\IEEEpubidadjcol

Calligraphic uppercase letters denote real \emph{Hilbert spaces}, as in $\mathcal{X}$, $\mathcal{V}$. 
We denote the \emph{scalar product} of a Hilbert space by $\langle \cdot , \cdot \rangle$ and the associated \emph{norm} by $\| \cdot \|$.
$2^\mathcal{V}$ denotes the \emph{power} set of $\mathcal{V}$, i.e., the set of all subsets of $\mathcal{V}$.
An \emph{operator} (or mapping) $\myop{A}: \mathcal{X} \to \mathcal{V}$ maps each point in $\mathcal{X}$ to a point in $\mathcal{V}$.
A \emph{set-valued} operator $\myop{A} : \mathcal{X} \to 2^\mathcal{V}$ maps each element in $\mathcal{X}$ to a set in $\mathcal{V}$.
$\textbf{I}$ denotes the \emph{identity} operator.
$\mathcal{B} (\mathcal{X}, \mathcal{V})$ denotes the space of \emph{bounded linear} operators from $\mathcal{X}$ to $\mathcal{V}$;
we set $\mathcal{B} (\mathcal{X}) \eqdef \mathcal{B} (\mathcal{X}, \mathcal{X})$. 
Given an operator $\myop{A} \in \mathcal{B} (\mathcal{X}, \mathcal{V})$, its \emph{adjoint} $\myop{A}^*$ is the operator $\myop{A}^*: \mathcal{V} \to \mathcal{X}$ such that for all $\myvec{x} \in \mathcal{X}$ and $\myvec{u} \in \mathcal{V}$, $\langle \myop{A} \myvec{x}, \myvec{u} \rangle = \langle \myvec{x}, \myop{A}^* \myvec{u} \rangle$.
$\mathcal{S} (\mathcal{X})$ denotes the space of self-adjoint bounded linear operators from $\mathcal{X}$ to $\mathcal{X}$, i.e., $\mathcal{S} (\mathcal{X}) \eqdef \{ \myop{A} \in \mathcal{B} (\mathcal{X}) \, | \, \myop{A} = \myop{A}^* \}$.
Given two operators $\myop{A}, \, \myop{B} \in \mathcal{S} (\mathcal{X})$, the Loewner partial ordering on $\mathcal{S} (\mathcal{X})$ is defined by $\myop{A} \succeq \myop{B} \Leftrightarrow \langle \myop{A} \myvec{x}, \myvec{x} \rangle \geq \langle \myop{B} \myvec{x}, \myvec{x} \rangle$, $\forall \, \myvec{x} \in \mathcal{X}$. 
An operator $\myop{A}$ is said to be \emph{positive semidefinite} if $\myop{A}$ is a self-adjoint bounded linear operator and $\myop{A} \succeq 0$. 
Let $\alpha \in [0, +\infty[$; $\mathcal{P}_\alpha (\mathcal{X})$ denotes the space of positive semidefinite operators $\myop{A}$ such that $\myop{A} \succeq \alpha \textbf{I}$, i.e., $\mathcal{P}_\alpha (\mathcal{X}) \eqdef \{ \myop{A} \in \mathcal{S}(\mathcal{X}) \, | \, \myop{A} \succeq \alpha \textbf{I} \}$.
Given an operator $\myop{A} \in \mathcal{P}_\alpha$, its \emph{positive square root} $\sqrt{\myop{A}}$ is the unique operator $\sqrt{\myop{A}} \in \mathcal{P}_\alpha$ such that $(\sqrt{\myop{A}})^2 = \myop{A}$.
For every $\myop{A} \in \mathcal{P}_\alpha(\mathcal{X})$, we define a \emph{semi}-scalar product and a \emph{semi}-norm (a scalar product and a norm if $\alpha > 0$) by $\langle \cdot, \cdot \rangle_\myop{A} \eqdef \langle \myop{A} \cdot , \cdot \rangle$ and by $\| \cdot \|_\myop{A} \eqdef \sqrt{\langle \myop{A} \cdot , \cdot \rangle}$, respectively. 
The \emph{domain} of a set-valued operator $\myop{A}: \mathcal{X} \to 2^{\mathcal{X}}$ is defined by $\text{dom } \myop{A} \eqdef \{ \myvec{x} \in \mathcal{X} \, | \, \myop{A} \myvec{x} \neq \emptyset \}$, its \emph{graph} by $\text{gra } \myop{A} \eqdef \{ (\myvec{x},\myvec{u}) \in \mathcal{X} \times \mathcal{X} \, | \, \myvec{u} \in \myop{A} \myvec{x} \}$, the set of \emph{zeros} by $\text{zer } \myop{A} \eqdef \{ \myvec{x} \in \mathcal{X} \, | \, \mathbf{0} \in \myop{A} \myvec{x} \}$, the \emph{range} of $\myop{A}$ by $\text{ran } \myop{A} \eqdef \{ \myvec{u} \in \mathcal{X} \, | \, (\exists \, \myvec{x} \in \mathcal{X}) \; \myvec{u} \in \myop{A} \myvec{x} \}$, and the \emph{inverse} of $\myop{A}$ by $\myop{A}^{-1} : \mathcal{X} \to 2^\mathcal{X} : \myvec{u} \to \{ \myvec{x} \in \mathcal{X} \, | \, \myvec{u} \in \myop{A} \myvec{x} \}$. 
We use the notation $\{\myvec{x}^k\}$ as a shorthand for representing the sequence $\{\myvec{x}^k\}_{k=1}^{+\infty}$.
We say that a sequence $\{\myvec{x}^k\}$ in $\mathcal{H}$ \emph{converges in the norm} (or \emph{strongly converges}) to a point $\myvec{x}$ in $\mathcal{H}$ if $\| \myvec{x}^k - \myvec{x} \| \to 0$ and say that it \emph{converges weakly} if, for every $\myvec{u} \in \mathcal{H}$,  $\langle \myvec{x}^k, \myvec{u} \rangle \to \langle \myvec{x}, \myvec{u} \rangle$. 
We denote weak convergence by $\xrightarrow{w}$. Strong convergence implies weak convergence to the same limit. In finite-dimensional spaces, weak convergence implies strong convergence.
The space of \emph{absolutely-summable sequences} in $\mathbb{R}$, i.e., the space of sequences $\{\myvec{x}^k\}$ in $\mathbb{R}$ such that $\sum_k |\myvec{x}^k| < \infty$, is denoted by $\spaceseq^1 (\mathbb{N})$; the set of summable sequences in $[0, + \infty [$ is denoted by $\spaceseq^1_+ (\mathbb{N})$.
We denote by $\bigoplus_{\itr \in \{1, \dots, \dimr \}}  \hilbtwo_\itr$ the Hilbert direct sum~\cite[Example 2.1]{Bauschke2011} of the Hilbert spaces $\hilbtwo_\itr$, $\itr \in \{1, \dots, \dimr \}$.
Given two set-valued operators $\myop{A} : \mathcal{X} \to 2^\mathcal{V}$ and $\myop{B} : \mathcal{X} \to 2^\mathcal{V}$, their \emph{parallel sum} is $\parsum{\myop{A}}{\myop{B}} \eqdef \left( \myop{A}^{-1} + \myop{B}^{-1} \right)^{-1}$. 
Additionally, we denote by $\mathbb{R}$ the set of real \emph{numbers}, by $\mathbb{R}^n$ the set of real column \emph{vectors} of length $n$, and by $\mathbb{R}^{m \times n}$ the set of real \emph{matrices} with $m$ rows and $n$ columns.
$\mathbf{a}^T$ denotes the transpose of a vector $\mathbf{a}$
and $\mathbf{A}^T$ denotes the transpose of a matrix $\mathbf{A}$.
$[\mathbf{a}]_i$ denotes the $i$-th element of a vector $\mathbf{a}$, $[\mathbf{A}]_{:j}$ denotes the $j$-th column of a matrix $\mathbf{A}$, and $[\mathbf{A}]_{ij}$ denotes the element in the $i$-th row and $j$-th column of a matrix $\mathbf{A}$.
$\| \mathbf{A} \|_F \; \eqdef \sqrt{\Tr(\mathbf{A} \mathbf{A}^*)}$ denotes the Frobenius norm of a matrix $\mathbf{A}$. 
Finally, let $f : \mathcal{X} \to ]-\infty,+\infty]$ be a \emph{function}. 
Its \emph{domain} is denoted by $\text{dom } f \eqdef \{ \myvec{x} \in \mathcal{X} \, | \, f(\myvec{x}) < + \infty \}$ and its \emph{epigraph} by $\text{epi } f \eqdef \{ (\myvec{x},s) \in \mathcal{X} \times \mathbb{R} \, | \, f(\myvec{x}) \leq s \}$.
The function $f$ is \emph{lower semi-continuous} if $\text{epi } f$ is \emph{closed} in $\mathcal{X} \times \mathbb{R}$, and \emph{convex} if $\text{epi } f$ is convex in $\mathcal{X} \times \mathbb{R}$. 
We use $\Gamma_0 (\mathcal{X})$ to denote the class of all lower semi-continuous convex functions $f$ from $\mathcal{X}$ to $]{-\infty},+\infty]$ that are \emph{proper}, i.e., such that $\text{dom } f \neq \emptyset$.
Given two functions $f \in \Gamma_0 (\mathcal{X})$ and $g \in \Gamma_0 (\mathcal{X})$, their \emph{infimal convolution} is $\infconv{f}{g}: \mathcal{X} \to [-\infty, +\infty] : \myvec{x} \to \underset{\myvec{u} \in \mathcal{X}}{\text{inf}} \left( f(\myvec{u}) + g(\myvec{x} - \myvec{u}) \right)$, where $\text{inf}$ denotes the infimum of a function. 
The \emph{indicator function} of a set $C \in \mathcal{X}$ is defined as $\delta_C(\myvec{x}) \eqdef 0$ if $\myvec{x} \in C$ and $\delta_C(\myvec{x}) \eqdef + \infty$ otherwise. 
When $C$ is convex, closed and non-empty, $\delta_C(\myvec{x}) \in \Gamma_0(\mathcal{X})$. 
The maximum and signum operators are denoted by $\text{max}(\cdot)$ and $\text{sgn}(\cdot)$, respectively.

\subsection{Background} \label{sec:background}

We assume that all the optimization problems under consideration are convex and have at least one minimizer. 
For the convenience of the reader, we collect some notions on convex analysis and monotone operators, and on semismooth Newton and active-set methods in Appendices~\ref{sec:cvx_analysis} and~\ref{sec:ssn}, respectively.
Consider now that we want to solve~\eqref{eq:split} with $N=0$, and that we use a splitting method for that effect, such as 
the \ac{FB} one. This method is characterized by the iteration
$\mathbf{x}^{k + 1} \leftarrow \text{prox}_{\tau^k g} \left(\mathbf{x}^k - \tau^k \nabla f \left(\mathbf{x}^k \right) \right)$, which is performed consecutively for all $k \in \mathbb{Z}$, and where $\tau^k > 0$.
The \ac{FB} method, also known as \emph{proximal--gradient}, is a first-order method and can be seen as an extension of the gradient method to nonsmooth optimization~\cite{Figueiredo2003, Daubechies2004}. 
It can be shown that the iterates $\mathbf{x}^{k}$ of this method converge to a solution of~\eqref{eq:split} at a sublinear rate, or, under certain assumptions, at a linear rate.
The inclusion of \emph{relaxation steps} in the \ac{FB} method may improve its convergence speed in practice; the resulting method is given by the iteration, for all $k \in \mathbb{Z}$, of
\begin{equation} \label{eq:iter_rfba}
	\mathbf{x}^{k + 1} \leftarrow \mathbf{x}^k + \lambda^k \left( \text{prox}_{\tau g} \left(\mathbf{x}^k - \tau \nabla f \left(\mathbf{x}^k \right)\right) - \mathbf{x}^k \right),
\end{equation}
where $\lambda^k > 0$ is the so-called \emph{relaxation parameter}~\cite{Combettes2005}. Additionally, it is also possible to consider the use of \emph{inertial steps}. An example of such use is the method characterized by the iteration
$\mathbf{x}^{k + 1} \leftarrow (1 - \alpha^k) \mathbf{x}^{k-1} + (\alpha^k- \lambda^k) \, \mathbf{x}^k + \lambda^k \text{prox}_{\tau^k g} \left[ \mathbf{x}^k - \tau^k \nabla f(\mathbf{x}^k) \right]$, $\forall \, k \in \mathbb{Z}$,
and with $\alpha^k > 0$. The convergence of this method was studied in~\cite{Bioucas2007} under some conditions on $f$ and $g$. The use of inertial steps is not limited to this specific formulation, and there are a number of alternatives~\cite{Daspremont2021}, most notably the one studied in~\cite{Beck2009b}, which was shown to obtain the optimal convergence rate, in function values $f (\mathbf{x}^k) + g (\mathbf{x}^k)$, for first-order methods.
The variations discussed so far have virtually the same computational cost as the original \ac{FB} method; this is not the case of the methods that we discuss next. The (relaxed) proximal--Newton method~\cite{Schmidt2011, Becker2012, Lee2014, TranDinh2015, Becker2019} makes use of second-order information about $f$. It is characterized by the iteration
$\mathbf{x}^{k + 1} \leftarrow \mathbf{x}^k + \lambda^k \left( \text{prox}_g^{\mathbf{B}^k} \left( \mathbf{x}^k - \left[ \mathbf{B}^k \right]^{-1} \nabla f \left( \mathbf{x}^k \right) \right) - \mathbf{x}^k \right)$,
where $\mathbf{B}^k$ is the Hessian of $f$ at $\mathbf{x}^k$. The local convergence rate of the iterates $\mathbf{x}^{k}$ generated by this method can be shown to be quadratic under some conditions. However, computing $\text{prox}_g^{\mathbf{B}^k}$ can be prohibitive in many problems of interest, and $\mathbf{B}^k$ is typically replaced by an approximation of the Hessian of $f$ at $\mathbf{x}^k$; it can be shown that the iterates of the resulting method converge to a solution at a superlinear rate under some conditions.
Finally, one can also consider incorporating a generalized derivative of $\text{prox}_{\tau^k g} \cdot$. This is the case of the methods known as semismooth, whose relation to~\eqref{eq:iter_rfba} we discuss in Section~\ref{sec:related}.

It is sometimes useful to consider splitting methods as instances of fixed-point methods. A solution $\hat{\mathbf{x}}$ of Problem~\eqref{eq:split} with $N=0$ satisfies the fixed-point equation $\hat{\mathbf{x}} = \text{prox}_{\tau g} \left(\hat{\mathbf{x}} - \tau \nabla f \left(\hat{\mathbf{x}} \right) \right)$, for $\tau > 0$, and it can be shown that the \ac{FB} method produces a sequence of points that is Fej\'{e}r monotone with respect to $\text{Fix } \text{prox}_{\tau g} [\cdot - \tau \nabla f(\cdot)]$~\cite{Combettes2009b}. 
More generally, if one wishes to find fixed points of a nonepansive operator $\mathbf{R}$, one can consider 
the {\ac{KM}} method, which is characterized by the iteration
$x^{k + 1} \leftarrow \mathbf{T}_{\lambda^k} \left( {x}^k \right) \eqdef {x}^k + \lambda^k (\mathbf{R} \left( {x}^k \right) - {x}^k)$,
where $\mathbf{T}_{\lambda^k}$ is termed a $\lambda^k$-\emph{averaged operator}.
By making $\mathbf{R} = \text{prox}_{\tau g} [\cdot - \tau \nabla f(\cdot)]$, we recover~\eqref{eq:iter_rfba}. 

\subsection{Contributions and outline} \label{sec:outline}

Consider yet another alternative to~\eqref{eq:iter_rfba}, where we replace the scalars $\lambda^k$ by linear operators $\ave^k$ such that, for every $\ite$, $\eye \succ \ave^\ite \succ 0$:
\begin{equation} \label{eq:ssn_intro}
	\mathbf{x}^{k + 1} \leftarrow \mathbf{x}^k + \ave^k \left( \text{prox}_{\tau^k g} \left[ \mathbf{x}^k - \tau^k \nabla f(\mathbf{x}^k) \right] - \mathbf{x}^k \right).
\end{equation}
The study of this iteration and, more generally, of the following extension of the \ac{KM} scheme:
\begin{equation} \label{eq:ekm}
	\varx^\iite \leftarrow \avop{\ave^\ite} \left( \varx^\ite \right) \eqdef \varx^\ite + \ave^\ite \left( \neop \left( \varx^\ite \right) - \varx^\ite \right),
\end{equation}
will be the basis of this work. For convenience, we call the operators $\avop{\ave^\ite}$, \emph{operator-averaged operators}. This formulation is rather general and allows one to consider extensions of current methods in two directions. Firstly, we consider a generic nonexpansive operator $\neop$ instead of just its instance as $\mathbf{R} = \text{prox}_{\tau g} [\cdot - \tau \nabla f(\cdot)]$, as is commonly done in the literature on semismooth Newton methods; this allows one to tackle more complex problems, such as ones of the form of~\eqref{eq:split} and Problem~\ref{eq:primalproblem2} [cf. Subsection~\ref{sec:pdprob}]. Secondly, we make very mild assumptions on the operators $\ave^\ite$ [cf. Section~\ref{sec:eKM}]; this allows one to consider a broad range of them, which may incorporate second-order information about the problem to solve (or not). We provide here what we consider to be an interesting example of this flexibility. Consider again that we want to solve~\eqref{eq:split} with $N=0$. We can adress this problem with an ``intermediate'' method between the proximal--Newton and the semismooth Newton methods if we make $\ave^\ite$ to be the inverse of the (possibly regularized) Hessian of $f$. This removes the difficulty of computing the operator $\text{prox}_g^{\mathbf{B}^k} \cdot$ and also of computing a generalized derivative of the typically nonconvex and nonsmoooth term $\text{prox}_{\tau g} \cdot$. This method is not, strictly speaking, a second-order method, but it seems to have fast convergence in practice [cf. Section~\ref{sec:apps}].

The outline of this work is as follows. In Section~\ref{sec:related}, we list work related to ours. In Section~\ref{sec:eKM}, we show that operator-averaged operators have a contractive property and prove the convergence of fixed-point iterations of these operators  under certain conditions. We base ourselves on the fact that they produce a sequence that is Fej\'{e}r monotone (specifically, variable-metric Fej\'{e}r monotone)~\cite{Combettes2013,Combettes2014b}. In more detail, we show in Subsection~\ref{sec:eVMFB} how operator-averaged operators can be used to extend an existing algorithm from the literature---a variable-metric \ac{FB} method, a generalization of the proximal--Newton discussed above---, and we prove its convergence under certain conditions. Additionally, we show in Subsection~\ref{sec:pdprob} how operator-averaged operators can be used to solve a primal--dual problem first studied in~\cite{Combettes2011b}, which generalizes many convex problems of interest in signal and image processing~\cite{Combettes2011b, Combettes2014b}. In Section~\ref{sec:apps}, we discuss applications of the proposed method to solve two inverse problems in signal processing, and discuss a framework that contemplates the incorporation of existing active-set methods, in a off-the-shelf fashion, to solve $\ell_2$-regularized minimization problems. 
We defer all proofs to Appendix~\ref{sec:proofs}. 

\section{Related work} \label{sec:related}

A similar formulation to~\eqref{eq:ekm} has been studied in the context of numerical analysis. Consider that one wishes to solve the system $\mathbf{A} \mathbf{x} = \mathbf{b}$, where $\mathbf{A} \in \mathbb{R}^{m \times n}$ and $\mathbf{b} \in \mathbb{R}^{n}$. This system can be solved for $\mathbf{x}$, under certain conditions, by iterating $\mathbf{x} ^\iite \leftarrow \mathbf{x} ^\ite + \lambda^\ite \left( \mathbf{A} \mathbf{x}^\ite - \mathbf{b} \right)$. This iteration, whose form is similar to the \ac{KM} scheme discussed above, is known as the Richardson iteration, and is a particular instance of a first-order linear nonstationary iterative method (examples of others are the Jacobi and the Gauss--Seidel iterations)~\cite{Axelsson1994,Ryabenkii2006}. Sometimes, for computational reasons (e.g., if $\mathbf{A}$ is deemed to be too ill conditioned), it is convenient to consider a preconditioned version of the linear system: $\mathbf{P} \mathbf{A} \mathbf{x} = \mathbf{P} \mathbf{b}$, where $\mathbf{P}$ is an invertible matrix. The corresponding preconditioned version of the Richardson iteration is $\mathbf{x} ^\iite \leftarrow \mathbf{x}^\ite + \lambda^\ite \mathbf{P} \left( \mathbf{A} \mathbf{x}^\ite - \mathbf{b} \right)$. 
Now consider that one wishes to solve the fixed-point equation of the operator $\neop$. One could equally consider a left-preconditioning scheme to solve this equation:
$\text{find } \varx \in \mathcal{X}$ such that $\ave \neop \left( \varx \right) = \ave \varx$.
By mimicking the preconditioned version of the Richardson iteration described above, one obtains the iteration 
$\varx^\iite \leftarrow \varx^\ite + \ave \left( \neop \left( \varx^\ite \right) - \varx^\ite \right)$,
which corresponds to making $\ave^\ite$ fixed in~\eqref{eq:ekm}, i.e., making, for all $\ite$, $\ave^\ite = \ave$, where $\ave \succ 0$.
The preconditioner $\ave$ is different from the operator $\mathbf{B}^k$ discussed in Subsection~\ref{sec:background}, and both can be used simultaneously: in Subsection~\ref{sec:eVMFB}, we consider a version of the \ac{FB} method that makes use of the two. The core iteration of that version of the method is
$\mathbf{x}^{k + 1} \leftarrow \mathbf{x}^k + \ave^k \left( \text{prox}_g^{\mathbf{B}^k} \left( \mathbf{x}^k - \left[ \mathbf{B}^k \right]^{-1} \nabla f \left( \mathbf{x}^k \right) \right) - \mathbf{x}^k \right)$.

Iterations of operator-averaged operators can also be seen as iterations of a line-search method if one considers $\lambda^\ite \left( \neop \left( \varx^\ite \right) - \varx^\ite \right)$ to be a step in the direction of the fixed-point residual $\neop \left( \varx^\ite \right) - \varx^\ite$ with step-length parameter $\lambda^\ite$. This idea has been explored in~\cite{Giselsson2016}, where the authors considered steps with length parameters $\lambda^\ite \geq 1$. The term $\ave^\ite \left( \neop \left( \varx^\ite \right) - \varx^\ite\right)$
can also be seen as indicating a search direction by noting the similarities with second-order line-search methods~\cite[Chapter 3]{Nocedal2006}: incorporating second-order information about the fixed-point residual in $\ave^\ite$ results in a Newton-like method. However, other directions could be taken by an appropriate selection of $\ave^\ite$. The idea of exploring directions different from the one given by $\neop \left( \varx^\ite \right) - \varx^\ite$ has also been explored in recent work~\cite{Themelis2019}. The authors studied an algorithm comprised of two steps: the first step corresponds to a step in a fast direction---making use of second-order information about the problem at hand---and the second one to a projection into a half-space. This projection is made to guarantee convergence, and enforces Fej\'{e}r monotonicity. The extension of the \ac{KM} scheme discussed in the present work could also be used to extend the method from that paper, since the use of operator-averaged operators expands the range of directions that still satisfy a Fej\'{e}r-monotonicity condition. 

Additionally, we can establish parallels between the proposed scheme~\eqref{eq:ekm} and coordinate-descent methods~\cite{Wright2015,Shi2016}. Consider the case where $\{\ave^\ite\}$ is a sequence of diagonal operators whose entries are either 0 or 1. These operators can then be used to select a particular coordinate---or block of coordinates---of $\varx$. Since the operator $\ave^\ite$ is binary, we have $\eye \succeq \ave^\ite \succeq 0$ and not $\eye \succ \ave^\ite \succ 0$, as we assumed initially. Algorithms resulting from this choice of binary $\ave^\ite$ have been studied elsewhere (e.g., \cite{Combettes2015,Latafat2019}) and are not analyzed in the remainder of this paper. 

We now discuss the connection with semismooth Newton and active-set methods alluded to previously.
If we make $\ave^\ite =  [\myop{V}(\mathbf{x}^k)]^{-1}$ and $\neop = \text{prox}_{\tau g} [\cdot - \tau \nabla f(\cdot)]$, where $f = \| \mathbf{y} - \mathbf{H} \cdot \|^2$, $g = \tau \| \cdot \|_1$, and $\myop{V}(\mathbf{x}^k)$ is defined to be the B-differential of $\neop$, we recover a semismooth Newton method to address the \ac{LASSO} problem [cf. Appendix~\ref{sec:ssn}]. 
Other active-set methods~\cite{Wen2010,Byrd2016,Chen2017} can also be considered as particular examples of the proposed scheme: in those cases, the operator $\ave^\ite$ takes a different form.

Results discussed in this paper were first presented in the PhD dissertation of the author~\cite[Subsections 5.2 and 5.4]{Simoes2017} where, to the best of our knowledge, the link between averaged operators and the class of methods comprising semismooth Newton and active-set ones was discussed for the first time. A preliminary version of the current work appeared in~\cite{Simoes2018}. Operator-averaged operators can be seen as instances of averaged-operators in non-Euclidean metrics [cf. Propostion~\ref{th:contractive}]; an analysis of methods to solve monotone inclusions based on averaged operators in variable metrics was independently proposed in~\cite{Glaudin2019}.

\section{An extension of the \ac{KM} scheme} \label{sec:eKM}

In what follows we consider exclusively finite-dimensional spaces; this is reasonable since we are mainly concerned with digital signals and images. We also focus on optimization, and we do not discuss a more general approach, often found in the literature, that makes use of monotone operators instead of gradients and subgradients of convex functions. However, in the proofs provided in Section~\ref{sec:proofs}, we avoided assumptions that prevented a generalization to infinite-dimensional spaces; we also developed the proofs in the setting of monotone operators.

\begin{definition}[Operator-averaged operators] \label{th:definition}
	Let $\nesubset$ be a nonempty subset of $\mathcal{X}$, let $\eps \in \; ]0, 1[$, and let $\ave$ be an operator in $\mathcal{X}$ such that
	$\Meig {\eye} \succeq \ave \succeq \meig {\eye}$, where $\Meig, \, \meig \in [\eps, 1-\eps]$.	
	We say that an operator $\avop{\ave}: \nesubset \to \mathcal{X}$ is an operator-averaged operator, or, more specifically, a $\ave$-averaged operator, if there exists a nonexpansive operator $\neop: \nesubset \to \mathcal{X}$ such that
	$\avop{\ave} \define ({\eye} - \ave) + \ave \neop$. 
\end{definition}

We have proved the following results:

\begin{proposition} \label{th:contractive}
	Let $\nesubset$ be a nonempty subset of $\mathcal{X}$, let $\eps \in \; ]0, 1[$, let $\ave$ be an operator in $\mathcal{X}$ satisfying
	$\Meig \emph{\eye} \succeq \ave \succeq \meig \emph{\eye}$, where $\Meig, \, \meig \in [\eps, 1-\eps]$.
	let $\neop: \nesubset \to \mathcal{X}$ be a nonexpansive operator, and let $\avop{\ave}: \nesubset \to \mathcal{X}$ be a $\ave$-averaged operator. Then the operator $\avop{\ave}$ is $\Meig$-averaged in the metric induced by $\inv{\ave}$. In other words, the operator $\avop{\ave}$ verifies
	$\norm{\avop{\ave} \left( \varx \right) - \avop{\ave} \left( \vary \right)}{\inv{\ave}}^2 \leq \norm{\varx - \vary}{\inv{\ave}}^2 - \frac{1 - \Meig}{\Meig} \norm{\left( \emph{\eye} - \avop{\ave} \right) \left( \varx \right) - \left( \emph{\eye} - \avop{\ave} \right) \left( \vary \right)}{\inv{\ave}}^2$
	for all $\varx, \, \vary \in \nesubset$.
\end{proposition}

\begin{theorem} \label{th:eKM}
	Let $\nesubset$ be a nonempty closed convex subset of $\mathcal{X}$, let $\eps \in \; ]0, 1[$, let $\seq{\seqave^\ite} \in \spaceseq^1_+(\mathbb{N})$, let $\seq{\ave^\ite}$ be a sequence of operators in $\mathcal{B} (\mathcal{X})$ such that, for all $\ite \in \mathbb{N}$,
	$\Meig^\ite \emph{\eye} \succeq \ave^\ite \succeq \meig^\ite \emph{\eye}$, with $\Meig^\ite, \, \meig^\ite \in [\eps, 1-\eps]$ and $\left( 1 + \seqave^\ite \right) \ave^\iite \succeq \ave^\ite$,	
	and let $\neop: \nesubset \to \nesubset$ be a nonexpansive operator such that $\emph{\fix} \neop \neq \emptyset$. Additionally, let $\varx^0 \in \nesubset$ and let, for all $\ite$, $\seq{\varx^\ite}$ be a sequence generated by~\eqref{eq:ekm}. Then $\seq{\varx^\ite}$ converges to a point in $\emph{\fix} \neop$.
\end{theorem}

\subsection{An extension of a variable-metric \ac{FB} method} \label{sec:eVMFB}

In this subsection, we show how operator-averaged operators can be used to extend proximal--Newton and variable-metric \ac{FB} methods of the type introduced in Subsection~\ref{sec:background}.
Consider Algorithm~\ref{algo:vmssn} to solve Problem~(\ref{eq:split}) with $N=0$. 
In what follows, $\seq{\fbvm^\ite}, \seq{\ave^\ite}$ are sequences of bounded linear operators, and $\seq{\errorout^\ite}, \seq{\errorin^\ite}$ are absolutely summable sequences that can be used to model errors.
\begin{algorithm} 
	\caption{} 	\label{algo:vmssn}
	\algsetup{indent=2em}	
	\begin{algorithmic}[1]
		\STATE Choose $\fbimp^0 \in \mathcal{X}$
		\STATE $\ite \leftarrow 1$
		\WHILE{stopping criterion is not satisfied}
		\STATE Choose $\fbscvm^\ite > 0,\ \fbvm^\ite \succ 0,\ \text{and } \eye \succ \ave^\ite \succ 0$
		\STATE $\fbexp^\ite \leftarrow \fbimp^\ite - \fbscvm^\ite \fbvm^\ite \left( \grad f(\fbimp^\ite) + \errorin^\ite \right)$ \label{eq:vmsnn_dual}
		\STATE $\fbimp^\iite \leftarrow \fbimp^\ite + \ave^\ite \left( \prox^{\inv{(\fbvm^\ite)}}_{\fbscvm^\ite g} \fbexp^\ite + \errorout^\ite - \fbimp^\ite \right)$ \label{eq:vmsnn_primal}
		\STATE $\ite \leftarrow \iite$
		\ENDWHILE
	\end{algorithmic}
\end{algorithm}
The following theorem establishes some convergence properties of this algorithm.

\begin{theorem} \label{th:ssn}
	Let $g \in \lsc (\hilb)$, let $\fbcoco \in \ ]0, + \infty[$, and let $f: \hilb \to \mathbb{R}$ be convex and differentiable with a $1/\fbcoco$-Lipschitzian gradient.		
Let $\seq{\fbvm^\ite}$ be a sequence of operators in $\hilb$ such that, for all $\ite \in \mathbb{N}$,
$\Meigvm \emph{\eye} \succeq \fbvm^\ite \succeq \meigvm \emph{\eye}$ and $\Meigvm, \, \meigvm \in \ ]0, + \infty[$,	
let $\fbeps \in \;]0, \min \{ 1, 2 \fbcoco / (\Meigvm + 1) \}[$, let $\seq{\ave^\ite}$ be a sequence of operators in $\hilb$ such that, for all $\ite$,
$\ave^\ite \fbvm^\ite = \fbvm^\ite \ave^\ite$, with $\Meig \emph{\eye} \succeq \ave^\ite \succeq \meig \emph{\eye}$ and $\Meig, \, \meig \in [\fbeps, 1]$,	
let $\seq{\fbseqvmave^\ite} \in \spaceseq^1_+(\mathbb{N})$, and let
$(1 + \fbseqvmave^\ite) \ave^\iite \fbvm^\iite \succeq \ave^\ite \fbvm^\ite$.	
Let $\seq{\fbscvm^\ite}$ be a sequence in $[\fbeps, (2 \fbcoco - \fbeps)/\Meigvm]$ and let $\seq{\errorout^\ite}, \seq{\errorin^\ite} \in \spaceseq^1(\mathbb{N})$. 
	
	Suppose that $\fbfixedset = \emph{\Argmin}(f + g) \neq \emptyset$. Let $\seq{\fbimp^\ite}$ be a sequence generated by Algorithm~\ref{algo:vmssn}. Then the following hold:	
	\begin{enumerate}[label=\emph{\arabic*)}]
		\item $\seq{\varx^\ite}$ is quasi-Fej\'{e}r monotone with respect to $\fbfixedset$ relative to $\seq{\inv{\left(\ave^\ite \fbvm^\ite\right)}}$. \label{th:fbfejer}
		\item $\seq{\varx^\ite}$ converges weakly to a point in $\fbfixedset$. \label{th:fbconv}
	\end{enumerate}
\end{theorem}

\begin{remark}
The assumption that, for every $\ite$, $\fbvm^\ite$ and $\ave^\ite$ commute, i.e., $\ave^\ite \fbvm^\ite - \fbvm^\ite \ave^\ite = 0$, may seem to be severe. However, take into account that existing algorithms consider one of these operators to be the identity operator: in semismooth Newton methods, $\fbvm^\ite = \eye$, $\forall \ite$, whereas in variable-metric \ac{FB} methods, $\ave^\ite = \eye$, $\forall \ite$. 
\end{remark}

\subsection{Primal--dual optimization algorithms} \label{sec:pdprob}

Combettes and Pesquet studied a primal--dual problem that generalizes many other problems~\cite[Problem 4.1]{Combettes2011b}. 
By being able to devise an algorithm to solve this problem, we are effectively tackling a large number of problems simultaneously---problem~\eqref{eq:split} is one of these.

Let $\cvxone \in \lsc(\hilb)$, let $\coco \in \; ]0, + \infty[$, let $\smooth: \hilb \to \mathcal{X}$ be convex and differentiable with a $\inv{\coco}$-Lipschitzian gradient, and let $\biasvar \in \hilb$. Let $\dimr$ be a strictly positive integer; for every $\itr \in \{ 1, \dots, \dimr \}$, let $\biasvarlnop_\itr \in \hilbtwo_\itr$, let $\cvxtwo_\itr \in \lsc(\hilbtwo_\itr)$, let $\strmaxmonopparam_\itr \in \; ]0, +\infty[$, let $\strong_\itr \in \lsc(\hilbtwo_\itr)$ be $\strmaxmonopparam_\itr$-strongly convex, let $\lnop_\itr \in \boundedlinop (\hilb, \hilbtwo_\itr)$ such that $\lnop_\itr \neq 0$, and let $\pdomega_\itr$ be real numbers in $]0, 1]$ such that $\sum_{\itr = 1}^{\dimr} \pdomega_\itr = 1$. The problem is as follows:

\begin{problem} \label{th:problem}

Solve the primal minimization problem,
\begin{equation*} \label{eq:primalproblem2} 
	\underset{\pr \in \mathcal{X}}{\text{{minimize}}} \, \cvxone (\pr) + \sum_{\itr=1}^{\dimr} \pdomega_\itr \left(\infconv{\cvxtwo_\itr}{\strong_\itr} \right) \left(\lnop_\itr \pr - \biasvarlnop_\itr \right) + \smooth (\pr) - \innerpro{\pr}{\biasvar}{},
\end{equation*}
together with its corresponding dual minimization problem,
\begin{align*} 
	\underset{\du_1 \in \mathcal{Y}_{m_1}, \cdots, \du_\itr \in \mathcal{Y}_{m_\itr}}{\text{{minimize}}} \, &\left( \infconv{\conj{\cvxone}}{\conj{\cvxtwo}} \right) \left( \biasvar - \sum_{\itr=1}^{\dimr} \pdomega_\itr  \conj{\lnop_\itr} \du_\itr \right) \\
	& \quad + \sum_{\itr=1}^{\dimr} \pdomega_\itr \left( \conj{\cvxtwo_\itr} (\du_\itr) + \conj{\strong_\itr} (\du_\itr) + \innerpro{\du_\itr}{\biasvarlnop_\itr}{} \right). \label{eq:dualproblem2}
\end{align*}	
The sets of solutions to these primal and dual problems are denoted by $P$ and $D$, respectively. We recover~\eqref{eq:split} by making $N=1$, $\biasvarlnop_1 = \myvec{0}$, $\biasvar = \myvec{0}$, and $\strong_1: \mathbf{u} \to 0$ if $\mathbf{u} = 0$ and $\strong_1: \mathbf{u} \to + \infty$ otherwise.

\end{problem}

Consider Algorithm~\ref{algo:stackevmfbapp} to solve Problem~\ref{th:problem}. In what follows, for all $\itr$, $\seq{\vm^\ite}$, $\seq{\ave^\ite}$, $\seq{\vm^\ite_\itr}$, $\seq{\ave^\ite_\itr}$ are sequences of linear operators, and $\seq{\errorpr^\ite}$, $\seq{\errordu^\ite_\itr}$, $\seq{\errorresolpr^\ite}$, $\seq{\errorresoldu^\ite_\itr}$ are absolutely-summable sequences that can be used to model errors. Algorithm~\ref{algo:stackevmfbapp} is an extension of~\cite[Example 6.4]{Combettes2014b}.
\begin{algorithm} 
\caption{} \label{algo:stackevmfbapp}
\algsetup{indent=2em}	
\begin{algorithmic}[1]
	\STATE Choose $\pr^0 \in \mathcal{X}$ and $\du_1^0 \in \mathcal{Y}_{m_1}, \cdots, \du_\itr^0 \in \mathcal{Y}_{m_\itr}$
	\STATE $\ite \leftarrow 1$
	\WHILE{stopping criterion is not satisfied}
	\FOR{$\itr = 1, \dots, \dimr$}
	\STATE Choose $\vm^\ite_\itr, \, \ave^\ite_\itr \succ 0 \text{ s.t. } \ave^\ite_\itr \prec \eye$
	\STATE $\pdu^\ite_\itr =\prox^{\inv{(\vm_\itr^\ite)}}_{\conj{\cvxtwo}_\itr} \big( \du^\ite_\itr + \vm_\itr^\ite \big( \lnop_\itr \pr^\ite$ \label{eq:stackevmfbapp_dual}
	\newline\makebox[2.5cm]{}$ - \grad{\conj{\strong}_\itr} \left(\du^\ite\right) - \errorresoldu^\ite_\itr - \biasvarlnop_\itr \big) \big) + \errordu^\ite_\itr$
	\STATE $\duo^\ite_\itr = 2 \pdu^\ite_\itr - \du^\ite_\itr$ \label{eq:stackevmfbapp_dual2}
	\STATE $\du^\iite_\itr = \du^\ite_\itr + \ave_\itr^\ite \left( \pdu^\ite_\itr - \du^\ite_\itr \right)$ \label{eq:ssn_in_algo2}
	\ENDFOR
	\STATE Choose $\vm^\ite, \, \ave^\ite \succ 0 \text{ s.t. } \ave^\ite \prec \eye$
	\STATE $\ppr^\ite = \prox_{\cvxone}^{\inv{(\vm^\ite)}} \big( \pr^\ite - \vm^\ite \big( \sum_{\itr = 1}^{\dimr} \pdomega_\itr \conj{\lnop}_\itr \duo^\ite_\itr$ \label{eq:stackevmfbapp_primal}
	\newline\makebox[3cm]{}$+ \grad{\smooth} \left(\pr^\ite\right) + \errorresolpr^\ite - \biasvar \big) \big) + \errorpr^\ite$
	\STATE $\pr^\iite =  \pr^\ite + \ave^\ite \left(\ppr^\ite - \pr^\ite \right)$ \label{eq:ssn_in_algo}
	\STATE $\ite \leftarrow \iite$
	\ENDWHILE
\end{algorithmic}	
\end{algorithm}
The following corollary establishes some convergence properties.

\begin{corollary} \label{th:stackevmfbapp}
	
Suppose that 
$\biasvar \in \emph{\ran} \left( \subgrad{\cvxone} + \sum_{\itr=1}^{\dimr} \pdomega_\itr \conj{\lnop_\itr} \left( \infconv{\subgrad{\cvxtwo}_\itr}{\subgrad{\strong}_\itr} \right) \left(\lnop_\itr \cdot - \biasvarlnop_\itr \right) + \grad{\smooth} \right)$
and set $\mincocomon \eqdef \min \{\coco, \strmaxmonopparam_1, \dots, \strmaxmonopparam_\dimr \}$. 
Let $\seq{\vm^\ite}$ be a sequence of operators in $\hilb$ such that, for all $\ite \in \mathbb{N}$,
$\Meigvm \emph{\eye} \succeq \vm^\ite \succeq \meigvm \emph{\eye}$ with $\Meigvm, \, \meigvm \in \ ]0, + \infty[$,	
let $\fbeps \in \;]0, \min \{ 1, \mincocomon \}[$, let $\seq{\ave^\ite}$ be a sequence of operators in $\hilb$ such that, for all $\ite$,
$\ave^\ite \vm^\ite = \vm^\ite \ave^\ite$ with $\Meig \emph{\eye} \succeq \ave^\ite \succeq \meig \emph{\eye}$ and  $\Meig, \, \meig \in [\fbeps, 1]$,	
and let
$\ave^\iite \vm^\iite \succeq \ave^\ite \vm^\ite$.	
Additionally, for every $\itr \in \{ 1, \dots, \dimr \}$, let $\seq{\vm^\ite_\itr}$ be a sequence of operators in $\hilbtwo_\itr$ such that, for all $\ite \in \mathbb{N}$,
$\Meigvm \emph{\eye} \succeq \vm^\ite_\itr \succeq \meigvm \emph{\eye}$ with $\Meigvm, \, \meigvm \in \ ]0, + \infty[$,	
let $\seq{\ave^\ite_\itr}$ be a sequence of operators in $\hilbtwo_\itr$ such that, for all $\ite$,
$\ave^\ite_\itr \vm^\ite_\itr = \vm^\ite_\itr \ave^\ite_\itr$ with $\Meig \emph{\eye} \succeq \ave^\ite_\itr \succeq \meig \emph{\eye}$ and $\Meig, \, \meig \in [\fbeps, 1]$,	
and let
$\ave^\iite_\itr \vm^\iite_\itr \succeq \ave^\ite_\itr \vm^\ite_\itr$.	
Let, for all $\itr$, $\seq{\errorpr^\ite}, \seq{\errordu^\ite}, \seq{\errorresolpr^\ite_\itr}, \seq{\errorresoldu^\ite_\itr} \in \spaceseq^1(\mathbb{N})$.	
For every $\ite$, set $\pddelta^\ite \define \isquareroot{\sum_{\itr = 1}^{\dimr} \pdomega_\itr  \norm{\squareroot{\vm^\ite_\itr} \lnop_\itr \squareroot{\vm^\ite}}{}^2} - 1$ and suppose that $\pdxi^\ite \define \frac{\pddelta^\ite}{(1+\pddelta^\ite) \Meigvm} \geq \frac{1}{2 \mincocomon - \fbeps}$.
Let $\seq{\pr^\ite}$ be a sequence generated by Algorithm~\ref{algo:stackevmfbapp}. Then $\pr^\ite$ converges to a point in $P$ and $\left(\du^\ite_1, \dots, \du^\ite_\dimr \right)$ converges to a point in $D$.
\end{corollary}

\section{Experiments} \label{sec:apps}

\subsection{Inverse integration} \label{sec:inv_int}

\begin{figure}[!t]
\begin{center}
	\includegraphics[scale=.4]{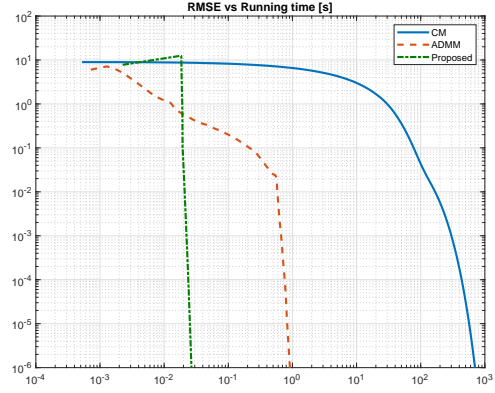}%
	\caption{RMSE, as a function of time, between the estimates of each iteration and the representative solution, for the three methods (\ac{CM}, \ac{ADMM}, and Proposed).}
	\label{fig:rmse_vs_time_inv_int}
\end{center}	
\end{figure}

We give a practical example of a simple problem that can be solved via Algorithm~\ref{algo:stackevmfbapp}. Consider the constrained problem
\begin{equation} \label{eq:lasso_constraint}
\underset{\mathbf{x} \in [c,d]^n}{\text{minimize}} \quad \| \mathbf{b} - \mathbf{H} \mathbf{x} \|^2 + \mu \| \mathbf{x} \|_1,
\end{equation}
where $\mathbf{b} \in \mathbb{R}^n$, $c \in \mathbb{R}$, $d \in \mathbb{R}$, $\mu>0$, $\mathbf{H} = {{}^{1}\!/_{n}} \widehat{\mathbf{H}}$, and $\widehat{\mathbf{H}} \in \mathbb{R}^{n \times n}$ is a lower-triangular matrix of ones. Griesse and Lorenz studied a non-constrained, and therefore simpler, version of this problem in the context of inverse integration~\cite[Section 4.1]{Griesse2008}. Problem~\eqref{eq:lasso_constraint} can be solved via Algorithm~\ref{algo:stackevmfbapp}---we discuss implementation details elsewhere~\cite[Section IV]{Simoes2018}.
We simulate an example similar to the one studied by Griesse and Lorenz~\cite[Section 4.1]{Griesse2008} but consider the noise to be Gaussian with a \ac{SNR} of 30 dB. We have set $\mu = 3 \times 10^{-3}$, $c=-80$, and $d=52$. We compared Algorithm~\ref{algo:stackevmfbapp} (denoted in what follows as \emph{Proposed}) with the \ac{ADMM} and with the \ac{CM}~\cite{Condat2013} to solve~\eqref{eq:lasso_constraint}. We manually tuned the different parameters of the three methods in order to achieve the fastest convergence results in practice. We arbitrarily chose the result of \ac{ADMM} after $10^7$ iterations as representative of the solution given by the three methods. Fig.~\ref{fig:rmse_vs_time_inv_int} illustrates the behavior of the three methods by showing the \ac{RMSE} 
between the estimates of each method and the representative solution, as a function of time. The three methods were initialized with the zero vector. The experiments were performed using MATLAB on an Intel Core i7 CPU running at 3.20~GHz, with 32~GB of RAM. 

In this example, we did not enforce any assumptions on $\ave^\ite$ but verified in practice that they were satisfied. However, in more complex examples, it may be necessary to devise a strategy that generates a sequence $\seq{\ave^\ite}$ satisfying the assumptions of Corollary~\ref{th:stackevmfbapp}. This is akin to the necessity of devising globalization strategies in other Newton-like methods~\cite[Chapter 8]{Facchinei2003b}.

\subsection{Plug-and-play methods}

Recently, a number of works have explored a technique known as \emph{plug-and-play \ac{ADMM}} to solve $\ell_1$-regularized problems with quadratic data-fitting terms. In that technique, the user substitutes one of the steps of the regular \ac{ADMM} by an off the shelf, state-of-the-art denoising algorithm~\cite{Chan2017, Romano2017}. Here, we explore a similar idea by considering the use of a semismooth Newton (or other active-set method) to replace a step of \ac{ADMM}. Whereas the concern of the former approach is to improve the \ac{SNR} of the estimates produced by it, our goal is instead to improve convergence speed. Consider problems of the form of~\eqref{eq:split} with $f=\frac{1}{2} \norm{\lnopH \pri - \mathbf{b}}{}^2$, $N=1$,
and where we further assume that $\cvxone (\pri) \in \Gamma_0(\mathbb{R}^n)$ is a sparsity-inducing regularizer.
Existing semismooth Newton methods do not consider the existence of a separate term $\cvxtwo$, and, in general, are not able to solve this problem efficiently, since that would require the computation of a (generalized) derivative of $\prox_{(g+h)}$. In what follows, we discuss an algorithm that can be used to solve this problem, that does not require that computation, and that allows the use of an existing semismooth Newton method in a plug-and-play fashion. 
We propose the following algorithm:
%
\begin{algorithm} 
\caption{} \label{algo:ssnal1}
\algsetup{indent=2em}	
\begin{algorithmic}[1]
	\STATE Choose $\pri^0 \in \mathbb{R}^n$, $\prispl^0 \in \mathbb{R}^n$, $\dua^0 \in \mathbb{R}^n$, $\tau > 0$
	\STATE $\ite \leftarrow 1$
	\WHILE{stopping criterion is not satisfied}
	\STATE $\ppri^\ite = \prox_{\tau \cvxone} \big( \pri^\ite - \tau \big( \grad{\smooth} (\pri^\ite)$ 
	\newline\makebox[3cm]{}$+ \gamma \left( \pri^\ite - \prispl^\ite + \dua^\ite \right) \big) \big)$ \label{algo:ssnal1ppri}
	\STATE $\pri^\iite =  \pri^\ite + \ave^\ite \left(\ppri^\ite - \pri^\ite \right)$ \label{algo:ssnal1pri}
	\STATE $\prispl^\iite = \prox_{\frac{\cvxtwo}{\gamma}} \left( \pri^\iite + \dua^\ite \right)$ \label{algo:ssnal1prspl}
	\STATE $\dua^\iite = \dua^\ite + \left( \pri^\iite - \prispl^\iite \right)$ \label{algo:ssnal1du}
	\STATE $\ite \leftarrow \iite$
	\ENDWHILE
\end{algorithmic}	
\end{algorithm}

Lines~\ref{algo:ssnal1ppri} and~\ref{algo:ssnal1pri} of Algorithm~\ref{algo:ssnal1} take a form similar to Algorithm~\ref{algo:vmssn} and can, in principle, be replaced by any active-set method~\cite{Wen2010,Muoi2013,Hans2015,Byrd2016,Chen2017,Stella2017,Themelis2019b}. 
In the following corollary, we discuss some of the convergence guarantees of this algorithm by showing that it is an instance of Algorithm~\ref{algo:stackevmfbapp}. 
Note that not all existing active-set methods will strictly obey the convergence conditions given here.

\begin{corollary} \label{th:plugnplayssn}
Suppose that 
$\myvec{0} \in \emph{\ran} \left( \subgrad{\cvxone} + \subgrad{\cvxtwo} + \grad{\smooth} \right)$.	
Set $\mincocomon =  \norm{\trans{\lnopH} \lnopH}{}$ and let $\fbeps \in \;]0, \min \{ 1, \mincocomon\}[$. 
For every $\ite \in \mathbb{N}$, set $\pddelta^\ite=\frac{1}{\tau \gamma} -1$  and suppose that 
$\frac{\pddelta^\ite}{(1+\pddelta^\ite) \Meigvm} \geq \frac{1}{2 \mincocomon - \fbeps}$ holds.
Let $\seq{\ave^\ite}$ be sequences of operators satisfying
$\Meig \emph{\eye} \succeq \ave^\ite \succeq \meig \emph{\eye}$,
$\Meig, \, \meig \in [\fbeps, 1]$,	
and let $\seq{\pri^\ite}$ be a sequence generated by Algorithm~\ref{algo:ssnal1}. Then $\pri^\ite$ converges weakly to a solution of~\eqref{eq:split} with $f=\frac{1}{2} \norm{\lnopH \pri - \mathbf{b}}{}^2$ and $N=1$.
\end{corollary}

\subsection{Spectral unmixing}

Hyperspectral images are images with a relatively large number of channels---usually known as \emph{spectral bands}---corresponding to short frequency ranges along the electromagnetic spectrum.
Frequently, their spatial resolution is low, and it is of interest to identify the materials that are present in a given pixel; a pixel typically corresponds to a mixture of different materials---known as \emph{endmembers}---, each with a corresponding spectral signature.
Spectral unmixing techniques produce a vector of \emph{abundances}, or percentages of occupation, for each endmember, in each pixel~\cite{Bioucas-Dias2012}.
Consider that we want to spectrally unmix a hyperspectral image with $N$ bands. We assume that the set of spectral signatures of the endmembers that may be present in a pixel is known through a database of $P$ materials (i.e., a database of reflectance profiles as a function of wavelength). We formulate the problem for a given pixel $j$ as
\begin{equation}
	\begin{aligned}
		& \underset{\mathbf{a}_j \in \mathbb{R}^{P}}{\text{minimize}} 
		& & \| [\mathbf{Y}_h]_{:j} - \mathbf{U} \mathbf{a}_j \|^2_2 + \mu \| \mathbf{a}_j \|_1\\
		& \text{subject to}
		& & [\mathbf{a}_j]_i \geq 0, \quad i = \{1, \cdots, P\}.
	\end{aligned}
\end{equation}
where $\mathbf{a}_j \in \R^{P}$ is the vector of each endmember's abundances for pixel $j$, to be estimated, $\mathbf{U} \in \R^{N \times P}$ is a matrix corresponding to the spectral database, $\mathbf{Y}_h \in \R^{N \times M}$ corresponds to a matrix representation of a hyperspectral image with $M$ pixels (i.e., corresponds to the lexicographical ordering of a 3-D data cube), and $\beta$ is a regularization parameter.
We tested two spectral dictionaries $\mathbf{U}$: a randomly generated dictionary and a real-world one. The former is given by sampling i.i.d. the standard Gaussian distribution, and the latter by a selection of 498 different mineral types from a USGS library, set up as detailed in~\cite{Bioucas2010}. The problem was generated as follows: we start by generating a vector of abundances with $P=224$ and with 5 nonzero elements, where the abundances are drawn from a Dirichlet distribution; we made $N=P$ and $M=100$, and added Gaussian noise such that it would result in a \ac{SNR} of 40 dB.
For this example, we implemented Algorithm~\ref{algo:ssnal1}, although we did not replace lines~\ref{algo:ssnal1ppri} and~\ref{algo:ssnal1pri} by any off-the-shelf method but by a direct implementation of those lines, which can be done as detailed in~\cite[Section IV]{Simoes2018}. Furthermore, we compared two versions of the proposed method, corresponding to different choices for the sequence of operators $\seq{\ave^\ite}$. The first, denoted as \emph{Proposed - variable}, corresponds to the inverse B-differential of $\left(\ppri^\ite - \pri^\ite \right)$; this choice is similar to the one made for the example detailed in Subsection~\ref{sec:inv_int}. The second, denoted as \emph{Proposed - fixed}, corresponds to making $\ave^\ite = \ave$ for all $\ite$, where $\ave$ corresponds to the scaled inverse of the regularized Hessian of $\| [\mathbf{Y}_h]_{:j} - \mathbf{U} \mathbf{a}_j \|^2_2$. By \emph{regularized}, we mean that we added the matrix $\epsilon \eye$ with $\epsilon = 10^2$ to the Hessian before computing the inverse, and by \emph{scaled}, that we scale the resulting inverse in order to guarantee that it obeys the condition $\Meig \eye \succeq \ave$; in other words, we made $\ave^{-1} = \rho \left( [\mathbf{Y}_h]_{:j}^T [\mathbf{Y}_h]_{:j} + \epsilon \eye \right)$, for a given $\rho$. We compared the two version of the proposed method as detailed in Subsection~\ref{sec:inv_int}.

\begin{figure}[tbhp]
\centering
\subfloat{
	\includegraphics[scale=.35,trim={40 20pt 40 0pt},clip]{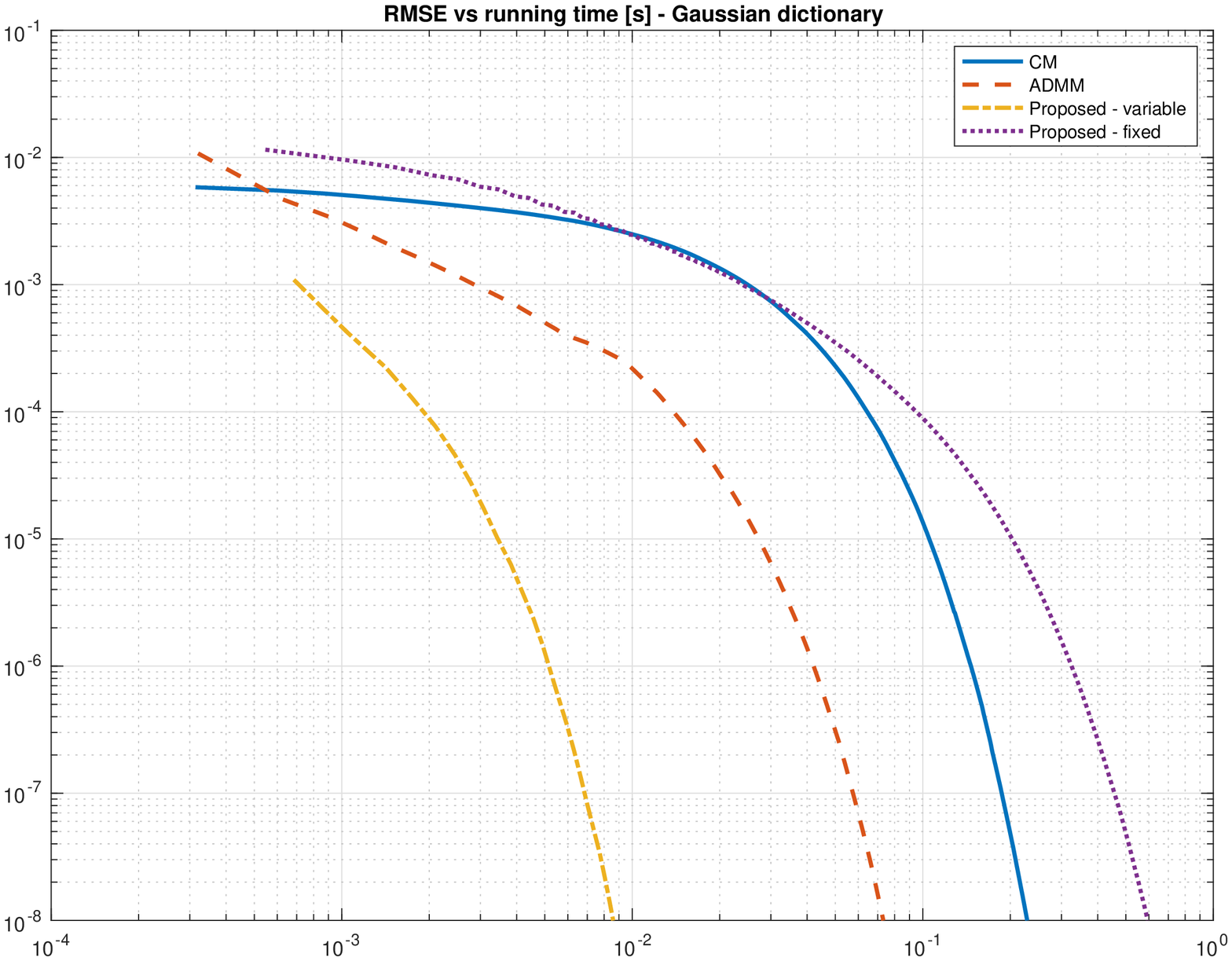}}

\subfloat{
	\includegraphics[scale=.35,trim={40 20pt 40 0pt},clip]{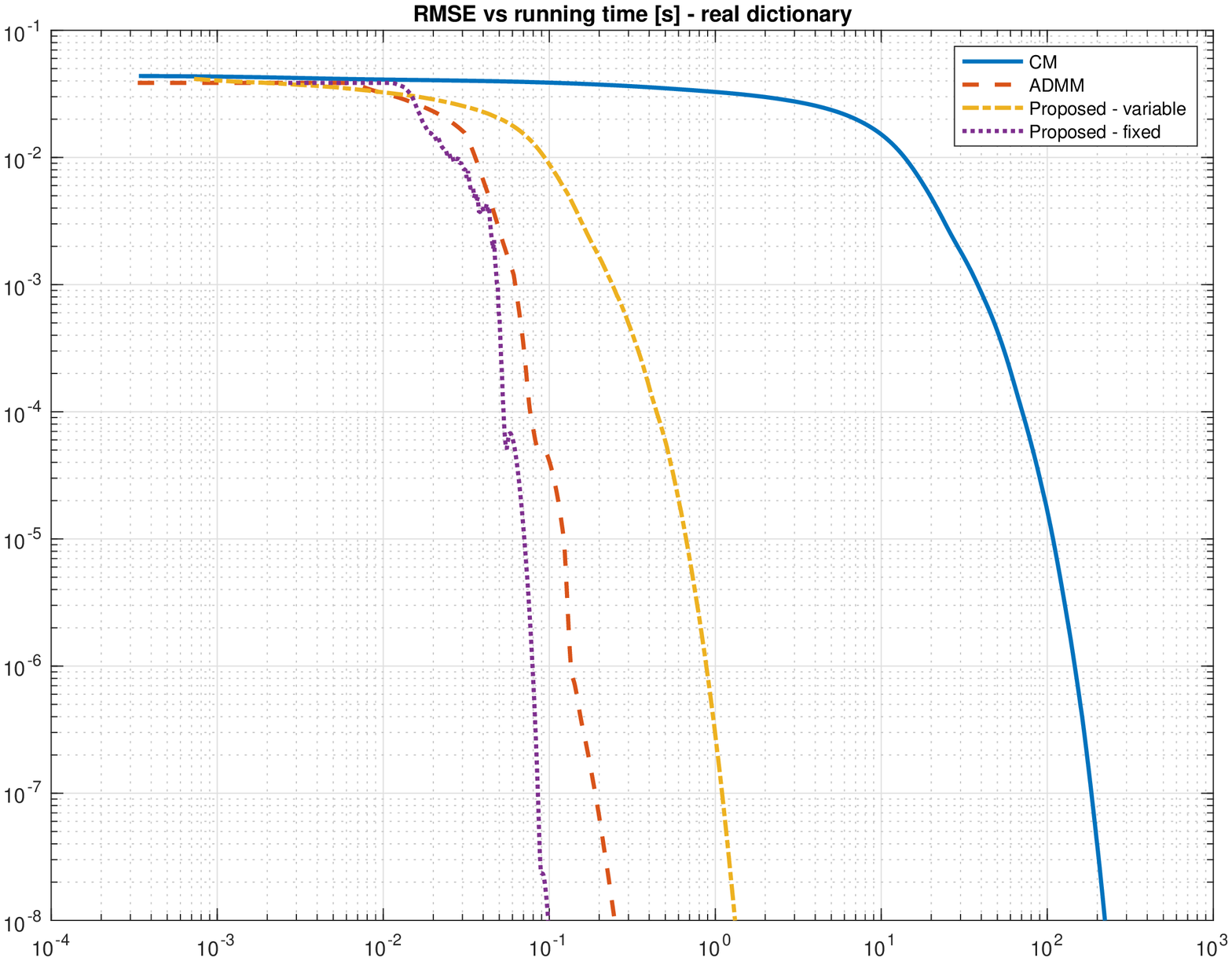}
}

\caption{RMSE, as a function of time, between the estimates of each iteration and the representative solution, for the four methods (\ac{CM}, \ac{ADMM}, Proposed - variable, and Proposed - fixed). Top: Gaussian dictionary $\mathbf{U}$; bottom: real dictionary $\mathbf{U}$.}
\label{fig:rmse_vs_time_spec_unmix}
\end{figure}

\section{Conclusions} \label{sec:conclusions}

This works discusses the use of operator-averaged operators to construct algorithms with fast convergence properties. These are particularly suitable to address problems with sparsity-inducing regularizers (and possibly other regularizers and constraints), and their behavior is similar to the one observed in the class of semismooth Newton and active-set methods if one selects an appropriate operator average. We tested two different choices of operators: one that corresponds precisely to the operators used in semismooth Newton methods and other that, while avoiding possibly hard-to-compute operations, still uses second-order information about the problem. The second choice corresponds to a method that is also different from the proximal--Newton one (and its corresponding approximations) and is, to be the best of our knowledge, new. We found experimentally that at least one of the two choices result in methods that are faster than conventional methods, but that they themselves may have significantly different convergence speeds. We hypothesize that the reason for this may be connected to the conditioning of the operator average---in the experiments in spectral unmixing, the matrices corresponding to the Gaussian and to the real dictionary have very different condition numbers, since the latter corresponds to the spectral profile of materials that are typically somewhat similar to each other. Future work directions to be explored are connected to the choice of an automatic strategy to select the sequence of $\ave^\ite$ operators, as well as possible ways to relax some of the assumptions on them.

\section*{Acknowledgments}
This work was supported by the Internal Funds KU Leuven project PDMT1/21/018.
I thank Prof.s Bioucas Dias, Borges de Almeida, and Panagiotis Patrinos for comments on how to improve this work. I dedicate this paper to the memory of Prof. Bioucas Dias.

\appendices

\section{On convex analysis and on monotone operators} \label{sec:cvx_analysis}

The notions of subgradient and subdifferential of a convex function (in the sense of Moreau and Rockafellar~\cite[Chapter 23]{Rockafellar1970}) are useful when dealing with nonsmooth functions. A vector $\mathbf{p} \in \mathbb{R}^n$ is said to be a \emph{subgradient} of a function $g \in \Gamma_0 (\mathbb{R}^n)$ at a point $\mathbf{x} \in \mathbb{R}^n$ if
$g(\mathbf{y}) \geq g(\mathbf{x}) + \langle \mathbf{p}, \mathbf{y} - \mathbf{x} \rangle$, $\forall \, \mathbf{y} \in \mathbb{R}^n$.
The set of all subgradients of $g$ at $\mathbf{x}$ is called the \emph{subdifferential of $g$ at $\mathbf{x}$} and is denoted by $\partial g (\mathbf{x})$. The set-valued operator $\partial g : \mathbb{R}^n \to 2^{\mathbb{R}^n} : \mathbf{x} \to \partial g (\mathbf{x})$ is called the \emph{subdifferential} of $g$. For a differentiable function $f \in \Gamma_0 (\mathbb{R}^n)$, the subdifferential at $\mathbf{x}$ is a singleton, i.e., $\partial f (\mathbf{x}) = \left\{ \nabla f (\mathbf{x}) \right\}$. The subdifferential operator is critical to our interests. 
We recall Fermat's rule~\cite[Theorem 16.2]{Bauschke2011}: $\mathbf{x}$ is a minimum of a proper convex function $g$ if and only if $\mathbf{0} \in \partial g (\mathbf{x})$. 
This means that the minimizers of $g$ are the zeros of the operator $\partial g$.
A smooth surrogate of a function $g \in \Gamma_0 (\mathbb{R}^n)$ is the so-called \emph{Moreau envelope}. It is defined by 
${}^\tau \! g (\mathbf{x}) \eqdef \inf_{\mathbf{u} \in \mathbb{R}^n} \left\{ g(\mathbf{u}) + \frac{1}{2 \tau} \| \mathbf{x} - \mathbf{u} \|^2 \right\}$,
where $\tau >0$. The function ${}^\tau \! g$ is continuously differentiable, even if $g$ is not. Both functions share the same minimizers~\cite[Proposition 12.9(iii)]{Bauschke2011}. The \emph{proximal operator} of $g$ is 
$\text{prox}_{\tau g}(\mathbf{x}) \eqdef \arg \min_{\mathbf{u} \in \mathbb{R}^n} \left\{ g(\mathbf{u}) + \frac{1}{2 \tau} \| \mathbf{x} - \mathbf{u} \|^2 \right\}$,
which is simply the point that achieves the infimum of the Moreau envelope (this point is unique, since $\mathbf{u} \to g(\mathbf{u}) + \frac{1}{2 \tau} \| \mathbf{x} - \mathbf{u} \|^2$ is strictly convex). The \emph{proximal operator} of $g$ \emph{relative to the norm} $\| \cdot \|^2_{\mathbf{U}}$, where $\mathbf{U} \in \mathbb{R}^{n \times n}$ and $\mathbf{U} \succ 0$, is 
$\text{prox}_{\tau g}^{\mathbf{U}} (\mathbf{x}) \eqdef \arg \min_{\mathbf{u} \in \mathbb{R}^n} \left\{ g(\mathbf{u}) + \frac{1}{2 \tau} \| \mathbf{x} - \mathbf{u} \|^2_{\mathbf{U}} \right\}$.
When $g (\mathbf{x}) = \| \mathbf{x} \|_1$, $\text{prox}_{\tau g} (\mathbf{x})$ can be evaluated component-wise. The proximal operator for each element reduces to the so-called \emph{soft-thresholding} operator~\cite{Donoho1995}, which is given by	$[\text{prox}_{\tau \| \cdot \|_1}(\mathbf{x})]_i = \text{max} \{| [\mathbf{x}]_i | - \tau, 0 \} \, \text{sgn} \left([\mathbf{x}]_i \right)$, $\forall i \in \{1, \cdots, n\}$.
The proximal operator
of the indicator function when the set $C \in \mathbb{R}^n$ is convex, closed and non-empty is an Euclidean projection onto $C$, which we denote by
$P_C (\mathbf{x}) \eqdef \arg \min_{\mathbf{u} \in C} \| \mathbf{x} - \mathbf{u} \|^2_2$.
For the proximal operators of more functions see, e.g., \cite{Combettes2011} or the website \url{http://proximity-operator.net/}.
Closely related to these ideas is the concept of Legendre--Fenchel \emph{conjugate} of a function $f$, which is defined as 
$f^* : \mathbb{R}^n \to [- \infty, + \infty] : \mathbf{x} \to \sup_{\mathbf{u} \in \mathbb{R}^n} \; \langle \mathbf{x} , \mathbf{u} \rangle - f(\mathbf{u})$.
We recall some of its properties. Consider that $f \in \Gamma_0(\mathbb{R}^n)$. Then $f^* \in \Gamma_0(\mathbb{R}^n)$ and, by the Fenchel--Moreau theorem~\cite[Theorem 13.32]{Bauschke2011}, the biconjugate of $f$ (the conjugate of the conjugate) is equal to $f$, i.e., $f^{**} = f$. Another property is that~\cite[Proposition 16.9]{Bauschke2011} $\mathbf{u} \in \partial f(\mathbf{x}) \Leftrightarrow \mathbf{x} \in \partial f^*(\mathbf{u})$, $\forall \, \mathbf{x}, \mathbf{u}  \in \mathbb{R}^n$,
or, in other words, $\partial f^* = (\partial f)^{-1}$. The notion of conjugate is also important in establishing the so-called \emph{Moreau's decomposition}, $\text{prox}_{\tau g} (\mathbf{x}) + \tau \, \text{prox}_{g^* / \tau} (\mathbf{x} / \tau) = \mathbf{x}$, $\forall \, \mathbf{x} \in \mathbb{R}^n$.
Finally, we say that $f$ is \emph{strongly convex} if $f - \frac{m}{2} \langle \mathbf{x}, \mathbf{x} \rangle$ is convex, with $m > 0$.

An operator $\myop{A}: \mathcal{X} \to 2^{\mathcal{X}}$ is said to be \emph{monotone} if
$\langle \myvec{u}-\myvec{v}, \myvec{x}-\myvec{y} \rangle \geq 0$, $\forall \, (\myvec{x},\myvec{u}) \in \text{gra } \myop{A}$, $\forall \, (\myvec{y},\myvec{v}) \in \text{gra } \myop{A}$.
An operator is \emph{maximally monotone} if there exists no other monotone operator whose graph properly contains $\text{gra } \myop{A}$. 
As an example, if $g \in \Gamma_0(\mathcal{X})$, than $\partial g$ is maximally monotone~\cite[Theorem 20.40]{Bauschke2011}. 
Let $\beta \in \, ]0, +\infty[$; we say that an operator is \emph{strongly monotone} with constant $\beta$ if $\myop{A} - \beta \eye$ is monotone. 
Monotone operators are connected to optimization problems as follows. 
Take, for example, Problem~(\ref{eq:split}) with $N=0$. 
According to Fermat's rule, its solutions should satisfy the inclusion $0 \in \nabla f (\myvec{x}) + \partial g (\myvec{x})$. 
Consequently, solving this problem can be seen as a particular case of the problem of finding a zero of the sum of two monotone operators $\myop{A}$ and $\myop{C}$ acting on a Hilbert space $\mathcal{X}$, i.e., 
$\text{find } x \in \mathcal{X}$, such that $0 \in \myop{A} (\myvec{x}) + \myop{C} (\myvec{x})$,
if ones makes $\myop{A}=\nabla f$ and $\myop{C}=\partial g$.
We now list some properties of operators. We say that an operator $\myop{A}$ is \emph{Lipschitz continuous} with constant $L > 0$ if 
$\| \myvec{u} - \myvec{v} \| \leq L \| \myvec{x} - \myvec{y} \|$, $\forall \, (\myvec{x},\myvec{u}) \in \text{gra } \myop{A}$, $\forall \, (\myvec{y},\myvec{v}) \in \text{gra } \myop{A}$.
When $L = 1$, the operator $\myop{A}$ is said to be \emph{nonexpansive}; when $L < 1$, it is said to be \emph{contractive}. 
Let $D$ be a nonempty set of $\mathcal{X}$ and let $\myop{R}: D \to \mathcal{X}$ be a nonexpansive operator. 
We say that an operator $\myop{A}: D \to \mathcal{X}$ is $\lambda$-\emph{averaged} if there exists $\lambda \in \; ]0,1[$ such that 
$\myop{A} = (1 - \lambda) \textbf{I} + \lambda \myop{R}$.
An averaged operator $\myop{A}$ obeys the following contractive property~\cite[Proposition 4.25]{Bauschke2011}:
$\norm{\myop{A} (\myvec{x}) - \myop{A} (\myvec{y})}{}^2 \leq \norm{\myvec{x} - \myvec{y}}{}^2 - \frac{1 - \lambda}{\lambda} \norm{\left( \textbf{I} - \myvec{A} \right) (\myvec{x}) - \left( \textbf{I} - \myop{A} \right) (\myvec{y})}{}^2$, $\forall \myvec{x} \in D$, $\forall \myvec{y} \in D$. 
When $\lambda = 1/2$, $\myop{A}$ is said to be \emph{firmly nonexpansive}. 
Proximal operators are examples of firmly nonexpansive operators~\cite[Corollary 23.8]{Bauschke2011}. 
Let $\beta \in \, ]0, +\infty[$; we say that a (single-valued) operator $\myop{C}: D \to \mathcal{X}$ is $\beta$-\emph{cocoercive} if
$\langle \myop{C}(\myvec{x}) - \myop{C}(\myvec{y}), \myvec{x} - \myvec{y} \rangle \geq \beta \|\myop{C}(\myvec{x}) - \myop{C}(\myvec{y})\|^2$, $\forall \, \myvec{x} \in D$, $\forall \, \myvec{y} \in D$.
An operator $\myop{C}$ is $\beta$-cocoercive if and only if $\beta \myop{C}$ is $\frac{1}{2}$-averaged~\cite[Remark 4.24(iv)]{Bauschke2011}. 
Let $f \in \Gamma_0(\mathcal{X})$ and let $\nabla f$ be $\beta$-Lipschitz continuous. 
Then, according to the Baillon--Haddad theorem~\cite[Corollary 18.16]{Bauschke2011}, $\nabla f$ is $\frac{1}{\beta}$-cocoercive.

In order to prove under which conditions iterative algorithms such as the ones that we discus in this work solve optimization problems, it can be useful to consider \emph{fixed-point} methods. 
The set of fixed points of an operator $\myop{A} : \mathcal{X} \to \mathcal{X}$ is
$\text{Fix } \myop{A} \eqdef \{ \myvec{x} \in \mathcal{X} \, | \, \myvec{x} = \myop{A} (\myop{x}) \}$.
If $\myop{A}$ is a Lipschitz-continuous operator, $\text{Fix } \myop{A}$ is closed~\cite[Proposition 4.14]{Bauschke2011}. 
If $\myop{A}$ is nonexpansive, $\text{Fix } \myop{A}$ is closed and convex~\cite[Corollary 4.15]{Bauschke2011}. 
Fixed-point methods try to find the fixed points of an operator (if they exist) by producing a sequence of points $\{\myvec{x}^{k}\}$ that should converge to one of them, given an initial point $\myvec{x}^0 \in \mathcal{X}$. 
A sequence $\{ \myvec{x}^{k}\}$ is said to be \emph{Fej\'{e}r monotone} with respect to a nonempty closed and convex set $S$ in $\mathcal{X}$ if
$\| \myvec{x}^{k+1} - \overline{\myvec{x}} \| \leq \| \myvec{x}^{k} - \overline{\myvec{x}} \|$, $\forall \, \overline{\myvec{x}} \in S$.
Such a sequence is bounded. 
Consequently, it possesses a subsequence that converges weakly to a point $\myvec{x} \in \mathcal{X}$. 
Such a point is said to be a \emph{weak sequential cluster point} of $\{\myvec{x}^{k}\}$, and we denote the set of weak sequential cluster points of $\{\myvec{x}^{k}\}$ by ${W}$. 
Interestingly, it is also a consequence of Fej\'{e}r monotonicity that a necessary and sufficient condition for the sequence $\{\myvec{x}^{k}\}$ to \emph{converge weakly} to a point in $S$ is that ${W} \subset S$~\cite{Combettes2009b},~\cite[Chapters 2 and 5]{Bauschke2011}. 
It is sometimes useful to consider the notions of quasi-Fej\'{e}r monotonicity and of Fej\'{e}r monotonicity relative to a variable metric. 
A sequence $\{\myvec{x}^{k}\}$ is said to be \emph{quasi}-Fej\'{e}r monotone with respect to a nonempty closed and convex set $S$ in $\mathcal{X}$ if
$\| \myvec{x}^{k+1} - \overline{\myvec{x}} \| \leq \| \myvec{x}^{k} - \overline{\myvec{x}} \| + \epsilon^k$,  $\exists \left\{\epsilon^\ite\right\} \in \ell^1_+(\mathbb{N})$, $\forall \, \overline{\myvec{x}} \in S$, 
and it said to be Fej\'{e}r monotone with respect to a nonempty closed and convex set $S$ in $\mathcal{X}$ \emph{relative to a sequence} $\left\{V^k\right\}$ 
if
$\| \myvec{x}^{k+1} - \overline{\myvec{x}} \|_{V^{k+1}} \leq (1 + \eta^k) \| \myvec{x}^{k} - \overline{\myvec{x}} \|_{V^k}$, $\forall \, \overline{\myvec{x}} \in S$
such that
(a) $V^k \in \mathcal{P}_\alpha(\mathcal{X})$ for $\alpha \in \, ]0, + \infty[$ and $\forall k \in \mathbb{N}$, (b) $\sup_k \| V^k \| < \infty$, and (c) $\left( 1 + \eta^k \right) V^{k+1} \succeq V^{k}$ with $\left\{\eta^\ite\right\} \in \ell^1_+(\mathbb{N})$ and $\forall k \in \mathbb{N}$.
The zeros of a monotone operator can be found by using fixed-point methods on appropriate operators. 
Let $\myop{A}: \mathcal{X} \to 2^{\mathcal{X}}$ be a maximally monotone operator and assume that $\text{zer } \myop{A} \neq \emptyset$.
Associated to this operator is its \emph{resolvent}
$J_{\tau \myop{A}} \eqdef (\textbf{I} + \tau \myop{A})^{-1}$.
The set of fixed points of $J_{\tau \myop{A}}$ coincides with the set of zeros of $\myop{A}$~\cite[Proposition 23.38]{Bauschke2011}. 
It can be shown that $J_{\tau \myop{A}}$ is firmly nonexpansive~\cite[Proposition 23.7]{Bauschke2011} and that $J_{\tau \partial g} = \text{prox}_{\tau g}$~\cite[Example 23.3]{Bauschke2011}. 
If one wishes to find a fixed point of a nonexpansive operator $\myop{R}$, one may use the \ac{KM} method [cf. Subsection~\ref{sec:background}].
Under certain conditions, it can be show that the sequnce $\{\myvec{x}^{k}\}$ is Fej\'{e}r monotone and that it converges weakly to a point in $\text{Fix } \myop{R}$, even when $\myop{R}$ is merely nonexpansive~\cite[Proposition 5.15]{Bauschke2011}. 
This method can be seen as a sequence of iterations of averaged operators since the operators $\myop{T}_{\lambda^k}$ are $\lambda^k$-averaged operators. 
In order to find the zeros of $\myop{A}$, one may consider a scheme based on an iteration of the form
$\myvec{x}^{k+1} = J_{\tau \myop{A}} \myvec{x}^k$, $\myvec{x}^0 \in \mathcal{X}$,
which converges since $J_{\tau A}$ is $1/2$-averaged. 

\section{On semismooth Newton methods} \label{sec:ssn}

We briefly discuss the use of semismooth Newton methods to solve the LASSO problem, which, we recall, consists in the minimization of the cost function
$\| \mathbf{H} \mathbf{x} - \mathbf{b} \|^2_2 + \mu \| \mathbf{x} \|_1$.
To address problems such as this, it was shown by Hinterm\"{u}ller~\cite{Hintermuller2003} that the use of certain semismooth Newton methods is equivalent to the use of \emph{active-set} ones (see below for a definition). 
The minimizer of the LASSO problem is assumed to be sparse; if we know which of its entries are zero, instead of solving this problem, we can solve an equivalent problem. 
Let $\mathbf{\hat{x}} \in \mathbb{R}^n$ be a solution to the problem, let $\hat{{A}}$ denote the set comprising the indices of the entries of $\mathbf{\hat{x}}$ that are zero, and let $\hat{{I}}$ denote the set comprising the remaining indices. 
The alternative problem 
\begin{equation} \label{eq:l1_ssn_active}
	\begin{aligned}
		& \underset{\mathbf{x} \in \mathbb{R}^{n}}{\text{minimize}}
		& & \| \mathbf{y} - \mathbf{H} \mathbf{x} \|^2_2 + \mu \| \mathbf{x} \|_1\\
		& \text{subject to}
		& & [\mathbf{x}]_i = 0, \, i \in \hat{{A}}
	\end{aligned}
\end{equation}
shares the same set of solutions. 
By writing the objective function of the latter problem as $\norm{\mathbf{y} - [\mathbf{H}]_{:j} [\mathbf{x}]_{j}}{2}^2 + \mu \norm{[\mathbf{x}]_{j}}{1}$ with $j \in \hat{{I}}$, it is clear that this problem has a much smaller number of non-constrained variables than the original. 
If operations involving $[\mathbf{H}]_{:j}$ are cheaper to perform than the ones involving the full matrix $\mathbf{H}$, an optimization algorithm will typically solve this problem faster than the original one. 
In practice, we do not know beforehand which entries of $\hat{\mathbf{x}}$ are zero. 
Active-set methods address this issue by finding estimates of the sets $\hat{{A}}$ and $\hat{{I}}$ by following some predefined strategy; 
the choice of strategy determining how both sets are estimated yields different algorithms~\cite{Wen2012,Milzarek2014,Solntsev2015,DeSantis2016,Liang2017,Li2018,Milzarek2019}. 
Define
$\myop{G}: \mathbb{R}^n \to \mathbb{R}^n: \mathbf{u} \to \mathbf{u} - \text{prox}_{\mu \tau \| \mathbf{u} \|_1} \left(\mathbf{u} - 2 \tau \mathbf{H}^* (\mathbf{H} \mathbf{u} - \mathbf{y}) \right)$.
The solution to the LASSO problem should satisfy the nonlinear equation $\myop{G} ({\mathbf{x}}) = \mathbf{0}$. 
This equation is nonsmooth, since $\text{prox}_{\mu \| \cdot \|_1}$ is not everywhere differentiable. 
There are, however, generalizations of the concept of differentiability that are applicable to nonconvex operators such as $\myop{G}$. 
One of them is the B(ouligand)-differential~\cite[Definition 4.6.2]{Facchinei2003a}, which is defined as follows. 
Suppose that a generic operator $\myop{G}: D \subset \mathbb{R}^n \to \mathbb{R}^m$ is locally Lipschitz, where $D$ is an open subset. 
Then by Rademacher's theorem, $\myop{G}$ is differentiable almost everywhere in $D$. 
Let $C$ denote the subset of $\mathbb{R}^n$ consisting of the points where $\myop{G}$ is differentiable (in the sense of Fr\'{e}chet~\cite[Definition 2.45]{Bauschke2011}). 
The \emph{B-differential of} $\myop{G}$ \emph{at} $\mathbf{x}$ is
$\partial_B \, \myop{G} (\mathbf{x}) \eqdef \left\{ \lim_{\mathbf{x}^j \to \mathbf{x}} \nabla \myop{G} \left(\mathbf{x}^j \right) \right\}$, 
where $\{\mathbf{x}^j\}$ is a sequence such that $\mathbf{x}^j \in C$ for all $j$ and $\nabla \myop{G} (\mathbf{x}^j)$ denotes the Jacobian of $\myop{G}$ at $\mathbf{x}^j$. 
As an example critical to our interests, consider the operator $\text{prox}_{\tau \| \cdot \|_1} (\mathbf{x})$. 
Its B-differential $\myop{B} = \partial_B \, \text{prox}_{\tau \| \cdot \|_1} (\mathbf{x})$ is given by~\cite[Proposition 3.3]{Griesse2008}
$[\myop{B}]_{ii} = 1$ if $\big| [\mathbf{x}]_i \big| > \tau$, $[\myop{B}]_{ii} = 0$ if $\big| [\mathbf{x}]_i \big| < \tau$, and $[\myop{B}]_{ii} \in \{0,1\}$ if $\big| [\mathbf{x}]_i \big| = \tau$.
Since the B-differential of an operator at a given point may not be unique, it may be convenient to consider a single $\mathbf{V} \in \myop{D}$, for example, the binary diagonal matrix
$[\mathbf{V}]_{ii} = 1$ if $\big| [\mathbf{x}]_i \big| > \tau$ and $[\mathbf{V}]_{ii} = 0$ otherwise. 

The generalization of the concept of differentiability just discussed can also be used to formulate the so-called semismooth Newton method based on the B-differential, which is characterized by the iteration 
$\mathbf{x}^{k+1} \leftarrow \mathbf{x}^{k} - [\mathbf{V}^k]^{-1} \, \myop{G}(\mathbf{x}^k)$,
where $\mathbf{V}^k \in \partial_B  \, G (\mathbf{x}^k)$. 
It can be shown that this method locally converges superlinearly for operators known as semismooth~\cite{Qi1993b}. 
Let $\mathbf{x} \in D$ and $\mathbf{d} \in \mathbb{R}^n$; semismooth operators are operators that are directionally differentiable at a neighborhood of $\mathbf{x}$ and that, for any $\mathbf{V} \in \partial_B \, \myop{G} (\mathbf{x}+\mathbf{d})$, satisfy the condition 
$\mathbf{V} \mathbf{d} - \myop{G}'(\mathbf{x}; \mathbf{d})=o(\|\mathbf{d}\|)$, $\mathbf{d} \to \mathbf{0}$,
where $\myop{G}'(\mathbf{x}; \mathbf{d})$ denotes the directional derivative~\cite[Definition 17.1]{Bauschke2011} of $\myop{G}$ at $\mathbf{x}$ along $\mathbf{d}$ and $o(\cdot)$ denotes little-O notation.
Examples of semismooth functions are the Euclidean norm and piecewise-differentiable functions~\cite[Chapter 2]{Ulbrich2011}, $\text{prox}_{\mu \| \cdot \|_1} (\mathbf{x})$ being an example of the latter.
For more details on these methods, see~\cite{Qi1999, Hintermuller2010, Ulbrich2009, Ulbrich2011}.

\section{Proofs} \label{sec:proofs}

This section includes the proofs of all the propositions, theorems, and corollaries of this work. It starts with a preliminary result.

\subsection*{Preliminary result} 

\begin{lemma} \label{th:pr1}
	Let $\alpha \in \; ]0, +\infty[$ and let $\iave \in \boundedlinopadjpos_\alpha(\hilb)$. Then, for all $\varx, \vary \in \hilb$,
	$\norm{\iave \varx + \left(\emph{\eye} - \iave \right) \vary}{}^2 = \innerpro{\iave \left(\iave - \emph{\eye} \right) \left( \varx - \vary \right)}{ \varx - \vary}{} + \norm{\varx}{\iave}^2 - \norm{\vary}{\iave}^2 + \norm{\vary}{}^2$.
\end{lemma}

\begin{proof}
	Fix $\varx$ and $\vary$ in $\hilb$. Then $\norm{\iave \varx + \left(\eye - \iave \right) \vary}{}^2 =$
	\begin{align*}
		&= \innerpro{\iave \varx}{\iave \varx}{} + 2 \innerpro{\iave \varx}{\left(\eye - \iave \right) \vary}{} \nonumber \\
		& \quad + \innerpro{\left(\eye - \iave \right) \vary}{\left(\eye - \iave \right) \vary}{} \nonumber \\
		&= \innerpro{\iave \varx}{\iave \varx}{} + 2 \innerpro{\iave \varx}{\vary}{} - 2 \innerpro{\iave \varx}{\iave \vary}{} + \innerpro{\vary}{\vary}{} \nonumber \\
		& \quad - 2 \innerpro{\iave \vary}{\vary}{} + \innerpro{\iave \vary}{\iave \vary}{} \nonumber \\
		&= \innerpro{\iave \varx}{\iave \varx}{} - 2 \innerpro{\iave \varx}{\iave \vary}{} + \innerpro{\iave \vary}{\iave \vary}{} + \innerpro{\iave \varx}{\varx}{} \nonumber \\ 		
		& \quad + \innerpro{\iave \vary}{\vary}{} - \innerpro{\iave \left(\varx - \vary \right)}{\varx - \vary}{} + \innerpro{\vary}{\vary}{} - 2 \innerpro{\iave \vary}{\vary}{} \nonumber \\
		&= \innerpro{\iave \left( \varx - \vary \right)}{\iave \left( \varx - \vary \right)}{} + \innerpro{\iave \varx}{\varx}{} - \innerpro{\iave \vary}{\vary}{} \nonumber \\
		& \quad - \innerpro{\iave \left( \varx - \vary \right)}{\varx - \vary}{} + \innerpro{\vary}{\vary}{} \nonumber \\
		&= \innerpro{\iave \left( \varx - \vary \right)}{\left(\iave - \eye \right)\left( \varx - \vary \right)}{} + \innerpro{\iave \varx}{\varx}{} - \innerpro{\iave \vary}{\vary}{}  \nonumber \\
		&\quad + \innerpro{\vary}{\vary}{}.
	\end{align*}
\end{proof}

\subsection*{Proof of Proposition~\ref{th:contractive}}

Fix $\varx$ and $\vary$ in $\nesubset$. By making $\iave \define \inv{\ave}$, we have 
$\inv{\meig} \eye \succeq \iave \succeq \inv{\Meig} \eye$.
and, by noting that $\neop = (\eye - \iave) + \iave \avop{\ave}$, we verify that $\norm{\neop \varx - \neop \vary}{}^2 = $
\begin{align}
	&= \norm{\left(\eye - \iave\right) \left(\varx - \vary\right) + \iave \left(\avop{\ave} \varx - \avop{\ave} \vary\right)}{}^2 \nonumber \\
	&\stackrel{(i)}{=} \big\langle \iave \left(\iave - \eye\right) \left( \left(\avop{\ave} - \eye \right) \varx - \left(\avop{\ave} - \eye \right) \vary \right), \nonumber \\
	& \qquad \left(\avop{\ave} - \eye \right) \varx - \left(\avop{\ave} - \eye \right) \vary \big\rangle \nonumber \\
	& \quad + \norm{\avop{\ave} \varx - \avop{\ave} \vary}{\iave}^2 - \norm{\varx - \vary}{\iave}^2 + \norm{\varx - \vary}{}^2 \nonumber \\
	&\stackrel{(ii)}{=} \norm{\left(\avop{\ave} - \eye \right) \varx - \left(\avop{\ave} - \eye \right) \vary }{\iave \left(\iave - \eye \right)}^2 \nonumber \\
	& \quad + \norm{\avop{\ave} \varx - \avop{\ave} \vary}{\iave}^2 - \norm{\varx - \vary}{\iave}^2 + \norm{\varx - \vary}{}^2,	
\end{align}
where step $(i)$ follows from Lemma~\ref{th:pr1} and step $(ii)$ follows from $\iave \left(\iave - \eye \right) \in \boundedlinopadjpos_0(\hilb)$, since $\inv{\meig} \geq \inv{\Meig} > 1$. The nonexpansiveness of $\neop$ implies that
$0 \geq \norm{\neop \varx - \neop \vary}{}^2 - \norm{\varx - \vary}{}^2 = \norm{\left(\avop{\ave} - \eye \right) \varx - \left(\avop{\ave} - \eye \right) \vary}{\iave \left(\iave - \eye \right)}^2 + \norm{\avop{\ave} \varx - \avop{\ave} \vary}{\iave}^2 - \norm{\varx - \vary}{\iave}^2$.
Consequently,
$\norm{\avop{\ave} \varx - \avop{\ave} \vary}{\iave}^2 \leq \norm{\varx - \vary}{\iave}^2 - \norm{\left(\avop{\ave} - \eye \right) \varx - \left(\avop{\ave} - \eye \right) \vary}{\iave \left(\iave - \eye \right)}^2$.
Since, for any given $\varz \in \hilb$, $\norm{\varz}{\iave \left(\iave - \eye \right)}^2 \geq \left( \inv{\Meig} - 1 \right) \norm{\varz}{\iave}^2$, then
$\norm{\avop{\ave} \varx - \avop{\ave} \vary}{\iave}^2 \leq \norm{\varx - \vary}{\iave}^2 - \left(\inv{\Meig} - 1\right) \norm{\left(\eye - \avop{\ave}\right) \varx - \left(\eye - \avop{\ave}\right) \vary}{\iave}^2$.
The claim follows by noting that $\left(\inv{\Meig} - 1\right) = \frac{1 - \Meig}{\Meig}$.

\subsection*{Proof of Theorem~\ref{th:eKM}}

\begin{enumerate}[label={\arabic*)}]
	\item Straightforward.
	
	\item Since $\varx^0 \in \nesubset$ and $\nesubset$ is convex, \eqref{eq:ekm} produces a well-defined sequence in $\nesubset$. By making $\iave^\ite \define \inv{(\ave^\ite)}, \forall \ite$, \cite[Lemma 2.1]{Combettes2013} yield, for all $\ite$, that
	$\inv{\left(\meig^\ite\right)} \eye \succeq \iave^\iite \succeq \inv{\left(\Meig^\ite\right)} \eye$ with
	$(1 + \seqave^\ite) \iave^\ite \succeq \iave^\iite$.	
	\eqref{eq:ekm} implies that, for all $\ite$, $\norm{\varx^\iite - \sol{\varx}}{\iave^\ite}^2 = \norm{\avop{\ave^\ite} \varx^\ite - \avop{\ave^\ite} \sol{\varx}}{\iave^\ite}^2 \stackrel{(i)}{\leq}$
	\begin{align}
		&\stackrel{(i)}{\leq} \norm{\varx^\ite - \sol{\varx}}{\iave^\ite}^2 \nonumber \\
		&\quad - \frac{1-\Meig^\ite}{\Meig^\ite} \norm{(\eye - \avop{\ave^\ite}) \varx^\ite + (\eye - \avop{\ave^\ite}) \sol{\varx}}{\iave^\ite}^2 \nonumber \\
		&= \norm{\varx^\ite - \sol{\varx}}{\iave^\ite}^2 - \frac{1-\Meig^\ite}{\Meig^\ite} \norm{\varx^\iite - \varx^\ite}{\iave^\ite}^2 \label{eq:avecontract} \\
		&\leq \norm{\varx^\ite - \sol{\varx}}{\iave^\ite}^2 \label{eq:avecontractfejer}, 
	\end{align}
	where step~$(i)$ follows from Proposition~\ref{th:contractive}. 
	Since, for any given $\varz \in \hilb$, we verify 
	$(1 + \seqave^\ite) \norm{\varz}{\iave^\ite}^2 \geq \norm{\varz}{\iave^\iite}^2$, $\forall \ite \in \mathbb{N}$,	
	using~\eqref{eq:avecontractfejer}, we get
	$\norm{\varx^\iite - \sol{\varx}}{\iave^\iite}^2 \leq (1 + \seqave^\ite) \norm{\varx^\iite - \sol{\varx}}{\iave^\ite}^2 \leq (1 + \seqave^\ite) \norm{\varx^\ite - \sol{\varx}}{\iave^\ite}^2$.
	
	\item Since $\seq{\varx^\ite}$ is Fej\'{e}r monotone with respect to $\fix \neop$ relative to $\seq{\iave^\ite}$, the sequence $\seq{\norm{\varx^\ite - \sol{\varx}}{\iave^\ite}^2}$ converges~\cite[Proposition 3.2(i)]{Combettes2013}. Define	
	$\fejersup \define \sup_\ite \norm{\varx^\ite - \sol{\varx}}{\iave^\ite} < + \infty$.		
	It follows from~\eqref{eq:ekm} that, for all $\ite$,
	\begin{align}
		\norm{\varx^\iite - \varx^\ite}{\iave^\ite}^2 &= \norm{\ave^\ite \left( \neop \varx^\ite - \varx^\ite \right)}{\iave^\ite}^2 \nonumber \\
		&=\norm{\left( \ave^\ite \right)^\frac{1}{2} \left( \neop \varx^\ite - \varx^\ite \right)}{}^2 \nonumber \\
		&\geq \meig^\ite \norm{\neop \varx^\ite - \varx^\ite}{}^2. \label{eq:resid_iter}
	\end{align}
	and, in view of~\eqref{eq:avecontract}, we can write 
	$\norm{\varx^\iite - \sol{\varx}}{\iave^\iite}^2 \leq$
	\begin{align}
		&\leq (1 + \seqave^\ite) \left( \norm{\varx^\ite - \sol{\varx}}{\iave^\ite}^2 - \frac{1-\Meig^\ite}{\Meig^\ite} \norm{\varx^\iite - \varx^\ite}{\iave^\ite}^2 \right) \nonumber \\
		&\leq (1 + \seqave^\ite) \norm{\varx^\ite - \sol{\varx}}{\iave^\ite}^2 - \left( 1-\Meig^\ite \right) \norm{\varx^\iite - \varx^\ite}{\iave^\ite}^2 \nonumber \\
		&\leq (1 + \seqave^\ite) \norm{\varx^\ite - \sol{\varx}}{\iave^\ite}^2 - \meig^\ite (1-\Meig^\ite) \norm{\neop \varx^\ite - \varx^\ite}{}^2 \nonumber \\
		&\leq \norm{\varx^\ite - \sol{\varx}}{\iave^\ite}^2 + \fejersup^2 \seqave^\ite - \eps^2 \norm{\neop \varx^\ite - \varx^\ite}{}^2. \label{eq:avecontractnonexp}
	\end{align}		
	
	For every $\Ite \in \mathbb{N}$, by iterating~\eqref{eq:avecontractnonexp} we can write that
$\eps^2 \sum_{\ite=0}^\Ite \norm{\neop \varx^\ite - \varx^\ite}{}^2 \leq \norm{\varx^0 - \sol{\varx}}{\iave^0}^2 - \norm{\varx^\Iite - \sol{\varx}}{\iave^\Iite}^2 + \sum_{\ite=0}^\Ite \fejersup^2 \seqave^\ite \leq \fejersup^2 + \sum_{\ite=0}^\Ite \fejersup^2 \seqave^\ite$.
	
	Since $\seq{\seqave^\ite}$ is absolutely summable, taking the limit as $\Ite \sconv + \infty$ yields
$\sum_{\ite=0}^{\infty} \norm{\neop \varx^\ite - \varx^\ite}{}^2 \leq \frac{1}{\eps^2} \left( \fejersup^2 + \sum_{\ite=0}^{\infty} \fejersup^2 \seqave^\ite \right) < \infty$.
	Consequently, $\neop \varx^\ite - \varx^\ite \sconv 0$.
	
	\item Let $\varx$ be a weak sequential cluster point of $\seq{\varx^\ite}$. It follows from~\cite[Corollary 4.18]{Bauschke2011} that $\varx \in \fix \neop$. In view of~\cite[Lemma 2.3]{Combettes2013} and~\cite[Theorem 3.3]{Combettes2013}, the proof is complete.
\end{enumerate}

\subsection*{Proof of Theorem~\ref{th:ssn}}

We start by proving a more general version of Algorithm~\ref{algo:vmssn}: consider that Lines~\ref{eq:vmsnn_dual} and~\ref{eq:vmsnn_primal} were replaced by $\fbexp^\ite \leftarrow \fbimp^\ite - \fbscvm^\ite \fbvm^\ite \left( \cocoercop \fbimp^\ite + \errorin^\ite \right)$ and $\fbimp^\iite \leftarrow \fbimp^\ite + \ave^\ite \left( \resol{\fbscvm^\ite \fbvm^\ite \maxmonresol} \fbexp^\ite + \errorout^\ite - \fbimp^\ite \right)$, respectively. Additionally, let $\maxmonresol: \hilb \to \powerset{\hilb}$ be a maximally monotone operator, let $\fbcoco \in \ ]0, + \infty[$, let $\cocoercop$ be a $\fbcoco$-cocoercive operator, and suppose instead that $\fbfixedset = \emph{\zer}(\maxmonop + \cocoercop) \neq \emptyset$. This more general version of Algorithm~\ref{algo:vmssn} allows one to address the problem of finding $\myvec{x} \in \mathcal{X}$ such that $\myvec{0} \in \myop{A} (\myvec{x}) + \myop{C} (\myvec{x})$.

Define, for all $\ite$,
$\maxmonresol^\ite \define \fbscvm^\ite \fbvm^\ite \maxmonresol$, 
$\cocoercop^\ite \define \fbscvm^\ite \fbvm^\ite \cocoercop$, 
$\fbvmave^\ite \define \fbvm^\ite \ave^\ite$, 
$\fbimpp^\ite \define \resol{\maxmonresol^\ite} \fbexp^\ite$,
$\fbimpq^\ite \define \resol{\maxmonresol^\ite} \left( \fbimp^\ite - \cocoercop^\ite \fbimp^\ite \right)$, 
$\fbimpwoe^\ite \define \fbimp^\ite + \ave^\ite \left( \fbimpq^\ite - \fbimp^\ite \right)$
We have from~\cite[Eq. (4.8)]{Combettes2014b} that
$\norm{\fbimpp^\ite - \fbimpq^\ite}{\inv{(\fbvm^\ite)}} \leq \frac{2 \fbcoco - \fbeps}{\sqrt{\Meigvm}} \norm{\errorin^\ite}{}$.
Additionally, for any $\fbfixed \in Z$, from~\cite[Eq. (4.12)]{Combettes2014b} we can write that
$\norm{\fbimpq^\ite - \fbfixed}{\inv{(\fbvm^\ite)}}^2 \leq \norm{\fbimp^\ite - \fbfixed}{\inv{(\fbvm^\ite)}}^2 - \fbeps^2 \norm{\cocoercop \fbimp^\ite - \cocoercop \fbfixed}{}^2 - \norm{(\fbimp^\ite - \fbimpq^\ite) - (\cocoercop^\ite \fbimp^\ite - \cocoercop^\ite \fbfixed)}{\inv{(\fbvm^\ite)}}^2$.

We now establish some identities. For all $\ite$, we have that
$\norm{\fbimp^\iite - \fbimpwoe^\ite}{\inv{(\fbvmave^\ite)}} = $
\begin{align} \label{eq:fbxsab}
	&= \norm{\left( \fbimp^\ite + \ave^\ite \left( \fbimpp^\ite + \errorout^\ite - \fbimp^\ite \right) \right) - \left(\fbimp^\ite + \ave^\ite \left( \fbimpq^\ite - \fbimp^\ite \right) \right)}{\inv{(\fbvmave^\ite)}} \nonumber \\
	&= \norm{\ave^\ite \left(\fbimpp^\ite + \errorout^\ite - \fbimpq^\ite \right)}{\inv{(\fbvmave^\ite)}} \nonumber \\
	&\leq \norm{\ave^\ite \errorout^\ite}{\inv{(\fbvmave^\ite)}} + \norm{\ave^\ite \left(\fbimpp^\ite - \fbimpq^\ite \right)}{\inv{(\fbvmave^\ite)}} \nonumber \\
	&= \norm{\squareroot{\inv{\left(\fbvm^\ite\right)}} \squareroot{\ave^\ite} \errorout^\ite}{} + \norm{\squareroot{\ave^\ite} \left(\fbimpp^\ite - \fbimpq^\ite \right)}{\inv{(\fbvm^\ite)}} \nonumber \\
	&\leq \norm{\squareroot{\inv{\left(\fbvm^\ite\right)}}}{} \norm{\squareroot{\ave^\ite}}{} \norm{\errorout^\ite}{} \nonumber \\
	&\quad + \norm{\squareroot{\ave^\ite}}{\inv{(\fbvm^\ite)}} \norm{\fbimpp^\ite - \fbimpq^\ite}{\inv{(\fbvm^\ite)}} \nonumber \\
	&\leq \squareroot{\norm{\inv{\left(\fbvm^\ite\right)}}{}} \squareroot{\norm{\ave^\ite}{}} \norm{\errorout^\ite}{} \\ 
	&\quad + \squareroot{\norm{\inv{\left(\fbvm^\ite\right)}}{}} \squareroot{\norm{\ave^\ite}{}} \norm{\fbimpp^\ite - \fbimpq^\ite}{\inv{(\fbvm^\ite)}} \nonumber \\		
	& \! \! \! \leq \sqrt{\Meig} \left( \frac{1}{\sqrt{\meigvm}} \norm{\errorout^\ite}{} + \frac{2 \fbcoco - \fbeps}{\sqrt{\meigvm \Meigvm}} \norm{\errorin^\ite}{} \right) \nonumber \\
	&\leq \frac{1}{\sqrt{\meigvm}} \norm{\errorout^\ite}{} + \frac{2 \fbcoco - \fbeps}{\sqrt{\meigvm \Meigvm}} \norm{\errorin^\ite}{}
\end{align}
and that
$\norm{\fbimpwoe^\ite - \fbfixed}{\inv{(\fbvmave^\ite)}}^2 =$
\begin{align} \label{eq:fbsz}
	&= \norm{\left( \fbimp^\ite + \ave^\ite \left( \fbimpq^\ite - \fbimp^\ite \right) \right) - \fbfixed}{\inv{(\fbvmave^\ite)}}^2 \nonumber \\
	&= \norm{\left( \eye - \ave^\ite \right) \left(\fbimp^\ite - \fbfixed \right) + \ave^\ite (\fbimpq^\ite - \fbfixed)}{\inv{(\fbvmave^\ite)}}^2 \nonumber \\
	&\stackrel{(i)}{=} \innerpro{\ave^\ite \left( \ave^\ite - \eye \right) \left( \fbimpq^\ite - \fbimp^\ite \right)}{\fbimpq^\ite - \fbimp^\ite}{\inv{(\fbvmave^\ite)}} \nonumber \\
	&\quad + \norm{\fbimpq^\ite - \fbfixed}{\inv{(\fbvmave^\ite)} \ave^\ite}^2 - \norm{\fbimp^\ite - \fbfixed}{\inv{(\fbvmave^\ite)} \ave^\ite}^2  \nonumber \\
	&\quad + \norm{\fbimp^\ite - \fbfixed}{\inv{(\fbvmave^\ite)}}^2 \nonumber \\
	&= - \norm{\fbimpq^\ite - \fbimp^\ite}{\inv{(\fbvm^\ite)} (\eye - \ave^\ite)}^2 + \norm{\fbimpq^\ite - \fbfixed}{\inv{(\fbvm^\ite)}}^2 \nonumber \\
	&\quad - \norm{\fbimp^\ite - \fbfixed}{\inv{(\fbvm^\ite)}}^2 + \norm{\fbimp^\ite - \fbfixed}{\inv{(\fbvmave^\ite)}}^2 \nonumber \\
	&\leq \norm{\fbimpq^\ite - \fbfixed}{\inv{(\fbvm^\ite)}}^2 - \norm{\fbimp^\ite - \fbfixed}{\inv{(\fbvm^\ite)}}^2 + \norm{\fbimp^\ite - \fbfixed}{\inv{(\fbvmave^\ite)}}^2 \nonumber \\
	&\stackrel{(ii)}{\leq} - \fbeps^2 \norm{\cocoercop \fbimp^\ite - \cocoercop \fbfixed}{}^2 + \norm{\fbimp^\ite - \fbfixed}{\inv{(\fbvmave^\ite)}}^2 \nonumber \\
	&\quad - \norm{(\fbimp^\ite - \fbimpq^\ite) - (\cocoercop^\ite \fbimp^\ite - \cocoercop^\ite \fbfixed)}{\inv{(\fbvm^\ite)}}^2
\end{align}
where step $(i)$ follows from Lemma~\ref{th:pr1} and step $(ii)$ from~\cite[Eq. (4.12)]{Combettes2014b}. 

Since, for any given $\varz \in \hilb$ and for all $\ite$, from~\cite[Lemma 2.1]{Combettes2013} we verify that
$(1 + \fbseqvmave^\ite) \inv{(\fbvmave^\ite)} \succeq \inv{(\fbvmave^\iite)}$ and  
$(1 + \fbseqvmave^\ite) \norm{\varz}{\inv{(\fbvmave^\ite)}}^2 \geq \norm{\varz}{\inv{(\fbvmave^\iite)}}^2$,
using inequation~\eqref{eq:fbsz}, we get
$\norm{\fbimpwoe^\ite - \fbfixed}{\inv{(\fbvmave^\iite)}}^2 \leq$
\begin{align}
	& \leq (1 + \fbseqvmave^\ite) \norm{\fbimp^\ite - \fbfixed}{\inv{(\fbvmave^\ite)}}^2 - \fbeps^2 \norm{\cocoercop^\ite \fbimp - \cocoercop \fbfixed}{}^2  \nonumber \\
	&\quad - \norm{(\fbimp^\ite - \fbimpq^\ite) - (\cocoercop^\ite \fbimp^\ite - \cocoercop^\ite \fbfixed)}{\inv{(\fbvm^\ite)}}^2 \label{eq:fbszxzcocoer} \\
	&\leq (1 + \fbseqvmave^\ite) \norm{\fbimp^\ite - \fbfixed}{\inv{(\fbvmave^\ite)}}^2 \label{eq:fbszxz} \\
	&\leq \fbseqsup^2 \norm{\fbimp^\ite - \fbfixed}{\inv{(\fbvmave^\ite)}}^2, \label{eq:fbszxzsup}
\end{align}
where $\fbseqsup \define \sup_\ite \sqrt{1 + \fbseqvmave^\ite}$.
We also define
$\fbeps^\ite \define \fbseqsup \left( \frac{1}{\sqrt{\meigvm}} \norm{\errorout^\ite}{} + \frac{2 \fbcoco - \fbeps}{\sqrt{\meigvm \Meigvm}} \norm{\errorin^\ite}{} \right)$.

Finally, these identities yield
$\norm{\fbimp^\iite - \fbimpwoe^\ite}{\inv{(\fbvmave^\iite)}}^2 \leq (1 + \fbseqvmave^\ite) \norm{\fbimp^\iite - \fbimpwoe^\ite}{\inv{(\fbvmave^\ite)}}^2 \leq (\fbeps^\ite)^2$.

\ref{th:fbfejer} We are now able to prove quasi-Fej\'{e}r monotonicity of $\seq{\fbimp^\ite}$:
$\norm{\fbimp^\iite - \fbfixed}{\inv{(\fbvmave^\iite)}} \leq$
\begin{align} \label{eq:fbxz}
	&\leq \norm{\fbimp^\iite - \fbimpwoe^\ite}{\inv{(\fbvmave^\iite)}} + \norm{\fbimpwoe^\ite - \fbfixed}{\inv{(\fbvmave^\iite)}} \nonumber \\
	&\leq \squareroot{1 + \fbseqvmave^\ite} \norm{\fbimp^\iite - \fbimpwoe^\ite}{\inv{(\fbvmave^\ite)}} + \squareroot{1 + \fbseqvmave^\ite} \norm{\fbimp^\ite - \fbfixed}{\inv{(\fbvmave^\ite)}} \nonumber \\
	&\leq \fbeps^\ite + \squareroot{1 + \fbseqvmave^\ite} \norm{\fbimp^\ite - \fbfixed}{\inv{(\fbvmave^\ite)}} \nonumber \\
	&\leq (1 + \fbseqvmave^\ite) \norm{\fbimp^\ite - \fbfixed}{\inv{(\fbvmave^\ite)}} + \fbeps^\ite.
\end{align}

Since $\seq{\errorout^\ite}$ and $\seq{\errorin^\ite}$ are absolutely summable, $\sum_\ite \fbeps^\ite < + \infty$. From the assumptions, in view of~\cite[Proposition 4.1(i)]{Combettes2013}, we conclude that $\seq{\fbimp^\ite}$ is quasi-Fej\'{e}r monotone with respect to $\fbfixedset$ relative to $\seq{\inv{(\fbvmave^\ite)}}$.

\ref{th:fbconv} As a consequence of~\ref{th:fbfejer} and~\cite[Proposition 4.1(ii)]{Combettes2013}, $\seq{\norm{\fbimp^\ite - \fbfixed}{\inv{(\fbvmave^\ite)}}}$ converges. We define
$\fbfejersup \define \sup_\ite \norm{\fbimp^\ite - \fbfixed}{\inv{(\fbvmave^\ite)}} < + \infty$.

Moreover,
$\norm{\fbimp^\iite - \fbfixed}{\inv{(\fbvmave^\iite)}}^2 =$
\begin{align}
	&= \norm{\fbimp^\iite - \fbimpwoe^\ite + \fbimpwoe^\ite - \fbfixed}{\inv{(\fbvmave^\iite)}}^2 \nonumber \\
	&\leq \norm{\fbimpwoe^\ite - \fbfixed}{\inv{(\fbvmave^\iite)}}^2 + \norm{\fbimp^\iite - \fbimpwoe^\ite}{\inv{(\fbvmave^\iite)}}^2 \nonumber \\
	&\quad + 2 \norm{\fbimpwoe^\ite - \fbfixed}{\inv{(\fbvmave^\iite)}} \norm{\fbimp^\iite - \fbimpwoe^\ite}{\inv{(\fbvmave^\iite)}} \nonumber \\
	&\leq (1 + \fbseqvmave^\ite) \norm{\fbimp^\ite - \fbfixed}{\inv{(\fbvmave^\ite)}}^2 - \fbeps^2 \norm{\cocoercop \fbimp^\ite - \cocoercop \fbfixed}{}^2 \nonumber \\ 
	&\quad - \norm{(\fbimp^\ite - \fbimpq^\ite) - (\cocoercop^\ite \fbimp^\ite - \cocoercop^\ite \fbfixed)}{\inv{(\fbvm^\ite)}}^2 + 2 \fbseqsup \fbfejersup \fbeps^\ite + (\fbeps^\ite)^2 \nonumber \\
	&\leq \norm{\fbimp^\ite - \fbfixed}{\inv{(\fbvmave^\ite)}}^2 - \fbeps^2 \norm{\cocoercop \fbimp^\ite - \cocoercop \fbfixed}{}^2 \nonumber \\ 
	&\quad - \norm{(\fbimp^\ite - \fbimpq^\ite) - (\cocoercop^\ite \fbimp^\ite - \cocoercop^\ite \fbfixed)}{\inv{(\fbvm^\ite)}}^2 \\
	&\quad + \fbfejersup^2 \fbseqvmave^\ite + 2 \fbseqsup \fbfejersup \fbeps^\ite + (\fbeps^\ite)^2. \label{eq:fbvm_iterate}	
\end{align}

For every $\Ite \in \mathbb{N}$, by iterating~\eqref{eq:fbvm_iterate}, we can write that
$\fbeps^2 \sum_{\ite=0}^\Ite \norm{\cocoercop \fbimp^\ite - \cocoercop \fbfixed}{}^2 \leq$
\begin{align}
	&\leq \norm{\fbimp^0 - \fbfixed}{\inv{(\fbvmave^0)}}^2 - \norm{\fbimp^{I+1} - \fbfixed}{\inv{(\fbvmave^{I+1})}}^2 \\
	&\quad + \sum_{\ite=0}^\Ite \left( \fbfejersup^2 \fbseqvmave^\ite + 2 \fbseqsup \fbfejersup \fbeps^\ite + (\fbeps^\ite)^2 \right) \nonumber \\
	&\leq \fbfejersup^2 + \sum_{\ite=0}^\Ite \left( \fbfejersup^2 \fbseqvmave^\ite + 2 \fbseqsup \fbfejersup \fbeps^\ite + (\fbeps^\ite)^2 \right).
\end{align}
Taking the limit from this inequality as $\Ite \to + \infty$ yields
$\sum_{\ite=0}^\infty \norm{\cocoercop \fbimp^\ite - \cocoercop \fbfixed}{}^2 \leq \frac{1}{\fbeps^2} \left( \fbfejersup^2 + \sum_{n} \left( \fbfejersup^2 \fbseqvmave^\ite + 2 \fbseqsup \fbfejersup \fbeps^\ite + (\fbeps^\ite)^2 \right) \right) < + \infty,$
since $\seq{\fbseqvmave^\ite}$ and $\seq{\fbeps^\ite}$ are absolutely summable. Following a similar reasoning, we can show that
$\sum_{\ite=0}^\infty \norm{(\fbimp^\ite - \fbimpq^\ite) - (\cocoercop^\ite \fbimp^\ite - \cocoercop^\ite \fbfixed)}{\inv{(\fbvm^\ite)}}^2 < + \infty$.

Let $\fbimp$ be a weak sequential cluster point of $\seq{\fbimp^\ite}$. In view of~\cite[Theorem 3.3]{Combettes2013}, it remains to be shown that $\fbimp \in \fbfixedset$. Since the sequences $\seq{\norm{\cocoercop \fbimp^\ite - \cocoercop \fbfixed}{}^2}$ and $\seq{\norm{(\fbimp^\ite - \fbimpq^\ite) - (\cocoercop^\ite \fbimp^\ite - \cocoercop^\ite \fbfixed)}{\inv{(\fbvm^\ite)}}^2}$ are absolutely summable, using the same arguments as in~\cite[Eqs. (4.26)-(4.31)]{Combettes2014b}, then $- \cocoercop \fbimp \in \maxmonop \fbimp$, which is equivalent to $\fbimp \in \fbfixedset$.


Finally, by making $\maxmonop = \subgrad{g}$ and $\cocoercop = \grad{f}$, we recover the original Algorithm~\ref{algo:vmssn}. By~\cite[Theorem 20.40]{Bauschke2011}, $\subgrad{g}$ is maximally monotone and, by~\cite[Corollary 18.16]{Bauschke2011}, $\grad{f}$ is $\fbcoco$-cocoercive. Additionally, $\Argmin (f+g) = \zer (\subgrad{f} + \grad{g})$, by~\cite[Corollary 26.3]{Bauschke2011}. 

\subsection*{Proof of Corollary~\ref{th:stackevmfbapp}}

The proof provided here follows the structure of similar proofs~\cite{Combettes2011b, Vu2011, Combettes2013}. As in the proof of Theorem~\ref{th:ssn}, we start by proving a more general version of Algorithm~\ref{algo:stackevmfbapp}: consider that Lines~\ref{eq:stackevmfbapp_dual} and~\ref{eq:stackevmfbapp_primal} were replaced by $\pdu^\ite_\itr \leftarrow \resol{\vm_\itr^\ite \inv{\maxmonopr}_\itr} \left( \du^\ite_\itr + \vm_\itr^\ite \left( \lnop_\itr \pr^\ite - \inv{\strmaxmonop}_\itr \du^\ite_\itr - \errorresoldu^\ite_\itr - \biasvarlnop_\itr \right) \right) + \errordu^\ite_\itr $ and $\ppr^\ite \leftarrow \resol{\vm^\ite \maxmonop} \left( \pr^\ite - \vm^\ite \left( \sum_{\itr = 1}^{\dimr} \pdomega_\itr \conj{\lnop}_\itr \duo^\ite_\itr + \cocoercop \pr^\ite + \errorresolpr^\ite - \biasvar \right) \right) + \errorpr^\ite$, respectively. 
Additionally, let $\maxmonop: \hilb \to \powerset{\hilb}$ be a maximally monotone operator, let $\coco \in \ ]0, + \infty[$, and let $\cocoercop: \hilb \to \hilb$ be $\coco$-cocoercive; for every $\itr \in \{ 1, \dots, \dimr \}$, let $\maxmonopr_\itr: \hilbtwo_\itr \to \powerset{\hilbtwo_\itr}$ be maximally monotone, let $\strmaxmonop_\itr: \hilbtwo_\itr \to \powerset{\hilbtwo_\itr}$ be maximally monotone and $\strmaxmonopparam_\itr$-strongly monotone, and suppose instead that 
$\biasvar \in {\ran} \left( \maxmonop + \sum_{\itr=1}^{\dimr} \pdomega_\itr \conj{\lnop_\itr} \left( \left( \parsum{\maxmonopr_\itr}{\strmaxmonop_\itr} \right) \left(\lnop_\itr \cdot - \biasvarlnop_\itr \right) \right) + \cocoercop \right)$.
This more general version of Algorithm~\ref{algo:stackevmfbapp} allows one to address the following primal--dual problem: solve the primal inclusion of finding ${\pr} \in \hilb$ such that $\biasvar \in \maxmonop {\pr} + \sum_{\itr=1}^{\dimr} \pdomega_\itr \conj{\lnop_\itr} \left( \left( \parsum{\maxmonopr_\itr}{\strmaxmonop_\itr} \right) \left(\lnop_\itr \pr - \biasvarlnop_\itr \right) \right) + \cocoercop {\pr}$ together with the dual inclusion of finding ${\du}_1 \in \hilbtwo_1, \dots, {\du}_\dimr \in \hilbtwo_\dimr$ such that $\exists \, \pr \in \hilb$ and $\biasvar - \sum_{\itr=1}^{\dimr} \pdomega_\itr \conj{\lnop_\itr} {\du}_\itr \in \maxmonop \pr + \cocoercop \pr$ and ${\du}_\itr \in \left( \parsum{\maxmonopr_\itr}{\strmaxmonop_\itr} \right) \left( \lnop_\itr \pr - \biasvarlnop_\itr \right)$, for all $\itr \in \{1, \dots, \dimr \}$. 

We start by introducing some notation. We denote by $\shilbtwo$ the Hilbert direct sum of the real Hilbert spaces $\hilbtwo_{\itr \in \{1, \dots, \dimr \}} $, i.e., $\shilbtwo = \bigoplus_{\itr \in \{1, \dots, \dimr \}}  \hilbtwo_\itr$. We endow this space with the following scalar product and norm, respectively:
$\innerpro{\cdot}{\cdot}{\shilbtwo} : (\stack{{a}}, \stack{{b}}) \mapsto 	\sum_{\itr=1}^{\dimr} \pdomega_\itr \innerpro{{a}_\itr}{{b}_\itr}{}$ and $\norm{\cdot}{\shilbtwo} : \stack{{a}} \mapsto \sqrt{\sum_{\itr=1}^{\dimr} \pdomega_\itr \norm{{a}_\itr}{}^2}$,
where $\stack{{a}} = ({a}_1, \dots, {a}_\itr), \stack{{b}} = ({b}_1, \dots, {b}_\itr) \in \shilbtwo$. Additionally, we denote by $\shilbthree$ the Hilbert direct sum $\shilbthree = \hilb \oplus \shilbtwo$ and endow the resulting space with the following scalar product and norm, respectively:
$\innerpro{\cdot}{\cdot}{\shilbthree} : \big( (\varx, \stack{\myvec{a}}), (\vary, \stack{\myvec{b}}) \big) \mapsto \innerpro{\varx}{\vary}{} + \innerpro{\stack{\myvec{a}}}{\stack{\myvec{b}}}{\shilbtwo}$ and $\norm{\cdot}{\shilbthree}: (\varx, \stack{\myvec{a}}) \mapsto \sqrt{\norm{\varx}{}^2 + \norm{\stack{\myvec{a}}}{\shilbtwo}}$,
where $\varx, \vary \in \hilb$.
We define, for all $\ite \in \mathbb{N}$,
$\sdur^\ite \in \shilbtwo \define \left(\du^\ite_1, \dots, \du^\ite_\dimr \right)$, $	
\spr^\ite \in \shilbthree \define \left(\pr^\ite, \sdur^\ite \right)$, $\sdu^\ite \in \shilbthree \define \left(\ppr^\ite, \pdu^\ite_1, \dots, \pdu^\ite_\dimr \right)$, $\serrorout^\ite \in \shilbthree \define \left(\errorpr^\ite, \errordu^\ite_1, \dots, \errordu^\ite_\dimr \right)$, $\serrorin^\ite \in \shilbthree \define \left(\errorresolpr^\ite, \errorresoldu^\ite_1, \dots, \errorresoldu^\ite_\dimr \right)$, $\svmerror^\ite \in \shilbthree \define \left(\inv{(\vm^\ite)} \errorpr^\ite, \inv{(\vm^\ite_1)} \errordu^\ite_1, \dots, \inv{(\vm^\ite_\dimr)} \errordu^\ite_\dimr \right)$.
We also define the operators
$\smaxmonresol : \shilbthree \to \powerset{\shilbthree} : (\varx, \stack{\myvec{a}}) \to \left( \sum_{\itr = 1}^{\dimr} \pdomega_\itr \conj{\lnop_\itr} a_\itr - \biasvar + \maxmonop \pr \right)$ 
$\times \left(- \lnop_1 \pr + \biasvarlnop_1 + \inv{\maxmonopr_1} a_1 \right)$ 
$\times \left(- \lnop_\dimr \pr + \biasvarlnop_\dimr + \inv{\maxmonopr_\dimr} a_\dimr \right)$, 
$\smaxmoncoco$ 
$: \shilbthree \to \shilbthree$
$: (\varx, \stack{\myvec{a}})$
$\to \left( \cocoercop \pr, \inv{\strmaxmonop}_1 a_1, \dots, \inv{\strmaxmonop}_\dimr a_\dimr \right)$, and
$\smaxmon : \shilbthree \to \shilbthree : (\varx, \stack{\myvec{a}}) \to \left( \sum_{\itr=1}^\dimr \pdomega_\itr \conj{\lnop_\itr} a_\itr, - \lnop_1 \pr, \dots, - \lnop_\dimr \pr \right)$,
and
$\svm^\ite : \shilbthree \to \shilbthree : (\varx, \stack{\myvec{a}}) \to \left( \vm^\ite \pr, \vm^\ite_1 a_1, \dots, \vm^\ite_\dimr a_\dimr \right)$, 
$\svmfb^\ite : \shilbthree \to \shilbthree : (\varx, \stack{\myvec{a}}) \to \Big( \inv{(\vm^\ite)} \pr + \sum_{\itr = 1}^{\dimr} \pdomega_\itr \conj{\lnop_\itr} a_\itr, \inv{(\vm^\ite_1)} a_1 + \lnop_1 \pr, \dots, \inv{(\vm^\ite_\dimr)} a_\dimr + \lnop_\dimr \pr \Big)$, 
$\sopt^\ite : \hilb \to \shilbtwo : \varx \to \left( \squareroot{\vm_1^\ite} \lnop_1 \pr, \dots, \squareroot{\vm_\dimr^\ite} \lnop_\dimr \pr \right)$, 
$\save^\ite : \shilbthree \to \shilbthree : (\varx, \stack{\myvec{a}}) \to \left( \ave^\ite \pr, \ave^\ite_1 a_1, \dots, \ave^\ite_\dimr a_\dimr \right)$,
$\svmfbave^\ite : \shilbthree \to \shilbthree : (\varx, \stack{\myvec{a}}) \to \left( \ave^\ite \vm^\ite \pr, \ave^\ite_1 \vm^\ite_1 a_1, \dots, \ave^\ite_\dimr \vm^\ite_\dimr a_\dimr \right)$,
where we note that the definition of the operator $\svmfb^\ite$ is not the same as the equivalent operator in~\cite[Eq. (6.10)]{Combettes2014b}.

We can further rewrite Lines~\ref{eq:stackevmfbapp_dual} and~\ref{eq:stackevmfbapp_dual2} of Algorithm~\ref{algo:stackevmfbapp} as $\inv{\left(\vm_\itr^\ite\right)} \left( \du^\ite_\itr - \pdu^\ite_\itr \right) + \lnop_\itr \pr^\ite - \inv{\strmaxmonop}_\itr \du^\ite_\itr \in \biasvarlnop_\itr + \inv{\maxmonopr}_\itr \left( \pdu^\ite_\itr - \errordu^\ite_\itr \right) + \errorresoldu^\ite_\itr - \inv{\left(\vm_\itr^\ite\right)} \errordu^\ite_\itr$ and Line~\ref{eq:stackevmfbapp_primal} as $\inv{\left(\vm^\ite\right)} \left(\pr^\ite - \ppr^\ite\right) + \sum_{\itr = 1}^{\dimr} \pdomega_\itr \conj{\lnop}_\itr \left( \du^\ite_\itr - \pdu^\ite_\itr \right) - \cocoercop \pr^\ite \in - \biasvar + \maxmonop \left( \ppr^\ite - \errorpr^\ite \right) + \sum_{\itr = 1}^{\dimr} \pdomega_\itr \conj{\lnop}_\itr \pdu^\ite_\itr + \errorresolpr^\ite - \inv{\left(\vm^\ite\right)} \errorpr^\ite$. 
In turn, these two modifications can be collapsed in the inclusion 
$\svmfb \left( \spr^\ite - \sdu^\ite \right) - \smaxmoncoco \spr^\ite \in \smaxmonresol (\sdu^\ite - \serrorout^\ite) + \smaxmon \serrorout^\ite + \serrorin^\ite - \svmerror^\ite$,
whereas Lines~\ref{eq:ssn_in_algo2} and~\ref{eq:ssn_in_algo} can be rewritten as
$\spr^\iite = \spr^\ite + \save^\ite \left( \sdu^\ite - \spr^\ite \right)$.

Set, for all $\ite$,	
$\serrorinpd^\ite \define \left(\smaxmon + \svmfb^\ite \right) \serrorout^\ite + \serrorin^\ite - \svmerror^\ite$.
Using the same arguments as in~\cite[Eqs. (3.25)-(3.30)]{Vu2011}, it can be shown that  
$\sdu^\ite = \resol{\inv{\left(\svmfb^\ite\right)} \smaxmonresol} \left( \spr^\ite - \inv{\left(\svmfb^\ite\right)} \left(\smaxmoncoco \spr^\ite + \serrorinpd^\ite \right) \right) + \serrorout^\ite$.
Note that the operator $\smaxmonresol$ is maximally monotone~\cite[Eqs. (3.7)-(3.9)]{Vu2011} and the operator $\smaxmoncoco$ is $\mincocomon$-cocoercive~\cite[Eq. (3.12)]{Vu2011}. Furthermore, it follows from the assumptions of the present corollary and~\cite[Lemma 3.1]{Combettes2014b} that
$\svm^\iite \succeq \svm^\ite \in \boundedlinopadjpos_\meigvm (\shilbthree)$ and $\norm{\inv{\left(\svm^\ite\right)}}{} \leq \frac{1}{\meigvm}$.
Since, for all $\ite$, $\svm^\ite \in \boundedlinopadj (\shilbthree)$, than $\svmfb^\ite \in \boundedlinopadj (\shilbthree)$. Additionally, by noting that $\smaxmon$ is linear and bounded, we verify, for all $\ite$, that
$\norm{\svmfb^\ite}{} \leq \norm{\inv{\left(\svm^\ite\right)}}{} + \norm{\smaxmon}{} \leq \frac{1}{\meigvm} + \sqrt{\sum_{\itr = 1}^{\dimr} \norm{\lnop_\itr}{}^2} \define \rho$
and that, for every $\pr \in \hilb$
$\norm{\sopt^\ite \pr}{\shilbtwo}^2 = \sum_{\itr = 1}^{\dimr} \pdomega_\itr \norm{\squareroot{\vm^\ite_\itr} \lnop_\itr \squareroot{\vm^\ite} \isquareroot{\vm^\ite} \pr}{}^2 \leq \norm{\pr}{\inv{(\vm^\ite)}}^2 \sum_{\itr = 1}^{\dimr} \pdomega_\itr \norm{\squareroot{\vm^\ite_\itr} \lnop_\itr \squareroot{\vm^\ite}}{}^2$ 
where $\mincocomon^\ite \define \sum_{\itr = 1}^{\dimr} \pdomega_\itr \norm{\squareroot{\vm^\ite_\itr} \lnop_\itr \squareroot{\vm^\ite}}{}^2, \forall_\ite$.
Then, following the arguments made in~\cite[Eq. 6.15]{Combettes2014b}, for every $\ite$ and every $\spr = (\varx,  {a}_1, \dots, {a}_\itr) \in \shilbthree$, we obtain $\innerpro{\svmfb^\ite \spr}{\spr}{\shilbthree} =$
\begin{align*}
	&= \norm{\pr}{\inv{(\vm^\ite)}}^2 + \sum_{\itr=1}^{\dimr} \pdomega_\itr \norm{a_\itr}{\inv{(\vm^\ite_\itr)}}^2 + 2 \Big\langle \isquareroot{(1+\pddelta^\ite) \mincocomon^\ite} \sopt^\ite \pr, \\
	& \qquad \squareroot{(1+\pddelta^\ite) \mincocomon^\ite} \left( \isquareroot{\vm^\ite_1} a_1, \dots, \isquareroot{\vm^\ite_\dimr} a_\dimr \right) \Big\rangle_{\shilbtwo} \\
	&\geq  \norm{\pr}{\inv{(\vm^\ite)}}^2 + \sum_{\itr=1}^{\dimr} \pdomega_\itr \norm{a_\itr}{\inv{(\vm^\ite_\itr)}}^2 \\
	&\quad - \left( \frac{\norm{\pr}{\inv{(\vm^\ite)}}^2}{1 + \pddelta^\ite} + (1 + \pddelta^\ite) \mincocomon^\ite \sum_{\itr=1}^{\dimr} \pdomega_\itr  \norm{a_\itr}{\inv{(\vm^\ite_\itr)}} \right) \\
	&\stackrel{(ii)}{=} \frac{\pddelta^\ite}{1+\pddelta^\ite} \left( \norm{\pr}{\inv{(\vm^\ite)}}^2 + \sum_{\itr=1}^{\dimr} \pdomega_\itr \norm{a_\itr}{\inv{(\vm^\ite_\itr)}}^2 \right) \\
	&\geq  \pdxi^\ite \norm{\pr}{\shilbthree}^2
\end{align*}
where step $(i)$ follows from the identity $2 \innerpro{\myvec{a}}{\myvec{b}}{} \geq - \norm{\myvec{a}}{}^2 - \norm{\myvec{b}}{}^2$ and step $(ii)$ follows from the fact that $(1 + \pddelta^\ite) \mincocomon^\ite = \frac{1}{1 + \pddelta^\ite}$. Following the arguments made in~\cite[Eqs. (6.16)-(6.18)]{Combettes2014b}, this last inequation implies that
$\sup_\ite \norm{\inv{\left(\svmfb^\ite\right)}}{} \leq 2 \mincocomon - \fbeps$ and $\inv{\left(\svmfb^\iite\right)} \succeq \inv{\left(\svmfb^\ite\right)} \in \boundedlinopadjpos_{1/ \rho} (\shilbthree)$,
It follows from the assumptions of the present corollary that
$\sup_\ite  \norm{\save^\ite}{} \leq 1$, $\save^\iite \succeq \save^\ite \in \boundedlinopadjpos_\meig (\shilbthree)$,
and
$\svmfbave^\iite \succeq \svmfbave^\ite$.
Moreover, it follows from~\cite[Lemma 3.1]{Combettes2014b} that $\sum_\ite \norm{\serrorout}{\shilbthree} \leq + \infty$, $\sum_\ite \norm{\serrorin}{\shilbthree} \leq + \infty$, and $\sum_\ite \norm{\svmerror}{\shilbthree} \leq + \infty$. 
It also follows that $\sum_\ite \norm{\serrorinpd}{\shilbthree} \leq + \infty$.
It is shown in~\cite[Eq. (3.13)]{Combettes2011b} that $\zer (\smaxmonresol + \smaxmoncoco) \neq \emptyset$. Additionally, following the arguments made in ~\cite[Eqs. (3.21) and (3.22)]{Combettes2011b}, if $(\sol{\pr}, \sol{\sdur}) \in \zer(\smaxmonresol + \smaxmoncoco)$, then $\sol{\pr} \in {P}$ and $\sol{\sdur} \in {D}$. 
Consequently, we have a $\sol{\spr} \define \left( \sol{\pr}, \sol{\du}_1, \dots, \sol{\du}_\itr \right)$ such that $\sol{\spr} \in \zer (\smaxmonresol + \smaxmoncoco)$ and $\spr^\ite \wconv \sol{\spr}$.	


By making, for every $\itr$, $\maxmonop = \subgrad{\cvxone}$, $\cocoercop = \grad{\coco}$, $\maxmonopr_\itr = \subgrad{\cvxtwo}_\itr$, and $\strmaxmonop_\itr=\subgrad{\strong_\itr}$, we recover the original Algorithm~\ref{algo:stackevmfbapp}. The current corollary is proven by using the same arguments as in~\cite[Theorem 4.2]{Combettes2011b}.

\subsection*{Proof of Corollary~\ref{th:plugnplayssn}}

Define, for every $\ite \in \mathbb{N}$,
$\duao^\ite = 2 \dua^\iite - \dua^\ite$,
$\bar{\dua}^\ite = \gamma \dua^\ite$,
$\bar{\duao}^\ite = \gamma \duao^\ite$,
and note that Lines~\ref{algo:ssnal1prspl} and~\ref{algo:ssnal1du} of Algorithm~\ref{algo:ssnal1} can be rewritten as 
$\dua^\iite = \dua^\ite + \pri^\iite - \prox_{\frac{\cvxtwo}{\gamma}} \left( \pri^\iite + \dua^\ite \right) = \frac{1}{\gamma} \prox_{\gamma \conj{\cvxtwo}} \left( \gamma \pri^\iite + \bar{\dua}^\ite \right).$
We can rewrite Lines~\ref{algo:ssnal1ppri}-\ref{algo:ssnal1du} of Algorithm~\ref{algo:ssnal1} by unfolding this algorithm:
$\prispl^\iite = \prox_{\frac{\cvxtwo}{\gamma}} \left( \pri^\ite + \dua^\ite \right)$, 
$\dua^\iite = \frac{1}{\gamma} \prox_{\gamma \conj{\cvxtwo}} \left( \gamma \pri^\ite + \bar{\dua}^\ite \right)$,
$\ppri^\ite = \prox_{\tau \cvxone} \left( \pri^\ite - \tau \left( \grad{\smooth} (\pri^\ite) + \gamma \left( \pri^\ite - \prispl^\iite + \dua^\iite \right) \right) \right)$, and
$\pri^\iite =  \pri^\ite + \ave^\ite \left(\ppri^\ite - \pri^\ite \right)$,
respectively.
Note that $\pri^\ite - \prispl^\iite + \dua^\iite = \pri^\ite - \prox_{\frac{\cvxtwo}{\gamma}} \left( \pri^\ite + \dua^\ite \right) + \frac{1}{\gamma} \prox_{\gamma \conj{\cvxtwo}_\itr} \left( \gamma \pri^\ite + \bar{\dua}^\ite \right) = \frac{2}{\gamma} \prox_{\gamma \conj{\cvxtwo}_\itr} \left( \gamma \pri^\ite + \bar{\dua}^\ite \right) - \dua^\ite = 2 \dua^\iite - \dua^\ite = \duao^\ite$.
This yields
$\ppri^\ite = \prox_{\tau \cvxone} \left( \pri^\ite - \tau \left( \gamma \duao^\ite + \grad{\smooth} (\pri^\ite) \right) \right)$.
We can likewise write that
$\bar{\duao}^\ite = 2\bar{\dua}^\iite - \bar{\dua}^\ite$
and
$\bar{\dua}^\iite = \prox_{\gamma \conj{\cvxtwo}} \left( \gamma \pri^\ite + \bar{\dua}^\ite \right)$.
By making $\dimr = 1,  \hilb = \mathbb{R}^n, \hilbtwo = \mathbb{R}^n, \lnop_1 = \mathbf{I}_n, \biasvarlnop_1 = \myvec{0}, \biasvar = \myvec{0}$, and $\forall_\ite, \vm^\ite_1 = \gamma \mathbf{I}_n, \vm^\ite = \tau \mathbf{I}_n, \errorresoldu^\ite_1 = \myvec{0}, \errordu^\ite_1 = \myvec{0}, \ave_{1}^\ite = \mathbf{I}_n, \errorresolpr^\ite = \myvec{0}, \errorpr^\ite = \myvec{0}$, it is clear that Algorithm~\ref{algo:ssnal1} is an instance of Algorithm~\ref{algo:stackevmfbapp}.
The current corollary is proven by invoking Corollary~\ref{th:stackevmfbapp}.

\bibliographystyle{IEEEtran}
\bibliography{refs}

\begin{thebibliography}{10}
\providecommand{\url}[1]{#1}
\csname url@samestyle\endcsname
\providecommand{\newblock}{\relax}
\providecommand{\bibinfo}[2]{#2}
\providecommand{\BIBentrySTDinterwordspacing}{\spaceskip=0pt\relax}
\providecommand{\BIBentryALTinterwordstretchfactor}{4}
\providecommand{\BIBentryALTinterwordspacing}{\spaceskip=\fontdimen2\font plus
\BIBentryALTinterwordstretchfactor\fontdimen3\font minus
  \fontdimen4\font\relax}
\providecommand{\BIBforeignlanguage}[2]{{%
\expandafter\ifx\csname l@#1\endcsname\relax
\typeout{** WARNING: IEEEtran.bst: No hyphenation pattern has been}%
\typeout{** loaded for the language `#1'. Using the pattern for}%
\typeout{** the default language instead.}%
\else
\language=\csname l@#1\endcsname
\fi
#2}}
\providecommand{\BIBdecl}{\relax}
\BIBdecl

\bibitem{Bauschke2011}
H.~Bauschke and P.~Combettes, \emph{Convex Analysis and Monotone Operator
  Theory in Hilbert Spaces}.\hskip 1em plus 0.5em minus 0.4em\relax New York,
  NY, USA: Springer, 2011.

\bibitem{Figueiredo2003}
M.~Figueiredo and R.~Nowak, ``An {EM} algorithm for wavelet-based image
  restoration,'' \emph{IEEE Trans. Image Process.}, vol.~12, no.~8, pp.
  906--916, Aug 2003.

\bibitem{Daubechies2004}
I.~Daubechies, M.~Defrise, and C.~De~Mol, ``An iterative thresholding algorithm
  for linear inverse problems with a sparsity constraint,'' \emph{Comm. Pure
  Appl. Math.}, vol.~57, no.~11, pp. 1413--1457, 2004.

\bibitem{Combettes2005}
P.~Combettes and V.~Wajs, ``Signal recovery by proximal forward--backward
  splitting,'' \emph{SIAM J. Multiscale Model. Simul.}, vol.~4, no.~4, pp.
  1168--1200, 2005.

\bibitem{Bioucas2007}
J.~M. Bioucas-Dias and M.~A.~T. Figueiredo, ``A new {TwIST}: Two-step iterative
  shrinkage/thresholding algorithms for image restoration,'' \emph{IEEE Trans.
  Image Process.}, vol.~16, no.~12, pp. 2992--3004, Dec 2007.

\bibitem{Daspremont2021}
\BIBentryALTinterwordspacing
A.~d'Aspremont, D.~Scieur, and A.~Taylor, ``{Acceleration Methods},'' in
  \emph{arXiv}, Mar. 2021. [Online]. Available:
  \url{https://arxiv.org/abs/2101.09545}
\BIBentrySTDinterwordspacing

\bibitem{Beck2009b}
A.~Beck and M.~Teboulle, ``A fast iterative shrinkage-thresholding algorithm
  for linear inverse problems,'' \emph{SIAM J. Imaging Sci.}, vol.~2, no.~1,
  pp. 183--202, 2009.

\bibitem{Schmidt2011}
M.~Schmidt, D.~Kim, and S.~Sra, ``Projected {Newton}-type methods in machine
  learning,'' in \emph{Optimization for Machine Learning}.\hskip 1em plus 0.5em
  minus 0.4em\relax MIT Press, 2011, pp. 305--330.

\bibitem{Becker2012}
S.~Becker and M.~Fadili, ``A quasi-{Newton} proximal splitting method,'' in
  \emph{Proc. 25th Int. Conf. Neural Informat. Process. Systems}, Lake Tahoe,
  Nevada, 2012, pp. 2618--2626.

\bibitem{Lee2014}
J.~Lee, Y.~Sun, and M.~Saunders, ``Proximal {Newton}-type methods for
  minimizing composite functions,'' \emph{SIAM J. Optim.}, vol.~24, no.~3, pp.
  1420--1443, 2014.

\bibitem{TranDinh2015}
Q.~Tran-Dinh, A.~Kyrillidis, and V.~Cevher, ``Composite self-concordant
  minimization,'' \emph{J. Mach. Learn. Res.}, vol.~16, no.~1, pp. 371--476,
  2015.

\bibitem{Becker2019}
S.~Becker, J.~Fadili, and P.~Ochs, ``On quasi-newton forward-backward
  splitting: Proximal calculus and convergence,'' \emph{SIAM J. Optim.},
  vol.~29, no.~4, pp. 2445--2481, 2019.

\bibitem{Combettes2009b}
P.~L. Combettes, ``Fej{\'e}r monotonicity in convex optimization,'' in
  \emph{Encyclopedia of Optimization}, C.~A. Floudas and P.~M. Pardalos,
  Eds.\hskip 1em plus 0.5em minus 0.4em\relax Boston, MA: Springer, 2009, pp.
  1016--1024.

\bibitem{Combettes2013}
P.~Combettes and B.~V\~{u}, ``Variable metric quasi-{F}ej\'{e}r monotonicity,''
  \emph{Nonlinear Anal.-Theor.}, vol.~78, pp. 17--31, 2013.

\bibitem{Combettes2014b}
P.~Combettes and B.~V{\~{u}}, ``Variable metric forward--backward splitting
  with applications to monotone inclusions in duality,'' \emph{Optim.},
  vol.~63, no.~9, pp. 1289--1318, 2014.

\bibitem{Combettes2011b}
P.~Combettes and J.-C. Pesquet, ``Primal--dual splitting algorithm for solving
  inclusions with mixtures of composite, {Lipschitzian}, and parallel-sum type
  monotone operators,'' \emph{Set-Valued Var. Anal.}, vol.~20, no.~2, pp.
  307--330, 2011.

\bibitem{Axelsson1994}
O.~Axelsson, \emph{Iterative Solution Methods}.\hskip 1em plus 0.5em minus
  0.4em\relax Cambridge University Press, 1994.

\bibitem{Ryabenkii2006}
V.~S. Ryaben'kii and S.~V. Tsynkov, \emph{A Theoretical Introduction to
  Numerical Analysis}.\hskip 1em plus 0.5em minus 0.4em\relax Boca Raton, FL,
  USA: CRC Press, 2006.

\bibitem{Giselsson2016}
P.~Giselsson, M.~F\"{a}lt, and S.~Boyd, ``Line search for averaged operator
  iteration,'' in \emph{IEEE Conf. Decision Control}, 2016, pp. 1015--1022.

\bibitem{Nocedal2006}
J.~Nocedal and S.~Wright, \emph{Numerical Optimization}.\hskip 1em plus 0.5em
  minus 0.4em\relax New York, NY, USA: Springer, 2006.

\bibitem{Themelis2019}
A.~Themelis and P.~Patrinos, ``Supermann: A superlinearly convergent algorithm
  for finding fixed points of nonexpansive operators,'' \emph{IEEE Trans.
  Automat. Contr.}, vol.~64, no.~12, pp. 4875--4890, 2019.

\bibitem{Wright2015}
S.~J. Wright, ``Coordinate descent algorithms,'' \emph{Math. Prog.}, vol. 151,
  no.~1, pp. 3--34, 2015.

\bibitem{Shi2016}
\BIBentryALTinterwordspacing
H.-J.~M. {Shi}, S.~{Tu}, Y.~{Xu}, and W.~{Yin}, ``{A Primer on Coordinate
  Descent Algorithms},'' in \emph{arXiv}, Jan. 2017. [Online]. Available:
  \url{https://arxiv.org/abs/1610.00040}
\BIBentrySTDinterwordspacing

\bibitem{Combettes2015}
P.~L. Combettes and J.-C. Pesquet, ``Stochastic quasi-{F}ej\'{e}r
  block-coordinate fixed point iterations with random sweeping,'' \emph{SIAM J.
  Optim.}, vol.~25, no.~2, pp. 1221--1248, 2015.

\bibitem{Latafat2019}
P.~Latafat, N.~M. Freris, and P.~Patrinos, ``A new randomized block-coordinate
  primal-dual proximal algorithm for distributed optimization,'' \emph{IEEE
  Trans. Automat. Contr.}, vol.~64, no.~10, pp. 4050--4065, 2019.

\bibitem{Wen2010}
Z.~Wen, W.~Yin, D.~Goldfarb, and Y.~Zhang, ``A fast algorithm for sparse
  reconstruction based on shrinkage, subspace optimization, and continuation,''
  \emph{SIAM J. Sci. Comput.}, vol.~32, pp. 1832--1857, 2010.

\bibitem{Byrd2016}
R.~H. Byrd, G.~M. Chin, J.~Nocedal, and F.~Oztoprak, ``A family of second-order
  methods for convex $\ell_1$-regularized optimization,'' \emph{Math.
  Program.}, vol. 159, no. 1-2, pp. 435--467, Sep. 2016.

\bibitem{Chen2017}
T.~Chen, F.~Curtis, and D.~Robinson, ``A reduced-space algorithm for minimizing
  $\ell_1$-regularized convex functions,'' \emph{SIAM J. Optim.}, vol.~27,
  no.~3, pp. 1583--1610, 2017.

\bibitem{Simoes2017}
\BIBentryALTinterwordspacing
M.~Sim\~{o}es, \emph{On some aspects of inverse problems in image
  processing}.\hskip 1em plus 0.5em minus 0.4em\relax Instituto Superior
  T\'{e}cnico, Universidade de Lisboa, Portugal \& Universit\'{e} Grenoble
  Alpes, France: PhD dissertation, 2017. [Online]. Available:
  \url{http://cascais.lx.it.pt/%7Emsimoes/dissertation/}
\BIBentrySTDinterwordspacing

\bibitem{Simoes2018}
M.~Simões, J.~Bioucas-Dias, and L.~B. Almeida, ``An extension of
  averaged-operator-based algorithms,'' in \emph{Europ. Sign. Proces. Conf.},
  Sep. 2018, pp. 752--756.

\bibitem{Glaudin2019}
L.~E. Glaudin, ``Variable metric algorithms driven by averaged operators,'' in
  \emph{Splitting Algorithms, Modern Operator Theory, and Applications}, H.~H.
  Bauschke, R.~S. Burachik, and D.~R. Luke, Eds.\hskip 1em plus 0.5em minus
  0.4em\relax Cham, Switzerland: Springer Int. Pub., 2019, pp. 227--242.

\bibitem{Griesse2008}
R.~Griesse and D.~Lorenz, ``A semismooth {N}ewton method for {T}ikhonov
  functionals with sparsity constraints,'' \emph{Inverse Probl.}, vol.~24,
  no.~3, p. 035007, 2008.

\bibitem{Condat2013}
L.~Condat, ``A primal--dual splitting method for convex optimization involving
  {L}ipschitzian, proximable and linear composite terms,'' \emph{J. Optim.
  Theory Appl.}, vol. 158, no.~2, pp. 460--479, 2013.

\bibitem{Facchinei2003b}
F.~Facchinei and J.-S. Pang, \emph{Finite-Dimensional Variational Inequalities
  and Complementarity Problems, Vols. II}.\hskip 1em plus 0.5em minus
  0.4em\relax Springer-Verlag, 2003.

\bibitem{Chan2017}
S.~H. Chan, X.~Wang, and O.~A. Elgendy, ``Plug-and-play {ADMM} for image
  restoration: Fixed-point convergence and applications,'' \emph{IEEE Trans.
  Comput. Imag.}, vol.~3, no.~1, pp. 84--98, March 2017.

\bibitem{Romano2017}
Y.~Romano, M.~Elad, and P.~Milanfar, ``The little engine that could:
  Regularization by denoising {(RED)},'' \emph{SIAM J. Imag. Sci.}, vol.~10,
  no.~4, pp. 1804--1844, 2017.

\bibitem{Muoi2013}
P.~Q. Muoi, D.~N. Hao, P.~Maass, and M.~Pidcock, ``Semismooth {Newton} and
  quasi-{Newton} methods in weighted $\ell_1$-regularization,'' \emph{J.
  Inverse Ill-Posed Probl.}, vol.~21, no.~5, pp. 665--693, 2013.

\bibitem{Hans2015}
E.~Hans and T.~Raasch, ``Global convergence of damped semismooth {N}ewton
  methods for {$\ell_1$} {T}ikhonov regularization,'' \emph{Inverse Probl.},
  vol.~31, no.~2, p. 025005, 2015.

\bibitem{Stella2017}
L.~Stella, A.~Themelis, and P.~Patrinos, ``Forward–backward quasi-{Newton}
  methods for nonsmooth optimization problems,'' \emph{Comput. Optim. Appl.},
  p. 443–487, 09 2017.

\bibitem{Themelis2019b}
A.~Themelis, M.~Ahookhosh, and P.~Patrinos, ``On the acceleration of
  forward-backward splitting via an inexact newton method,'' in \emph{Splitting
  Algorithms, Modern Operator Theory, and Applications}, H.~H. Bauschke, R.~S.
  Burachik, and D.~R. Luke, Eds.\hskip 1em plus 0.5em minus 0.4em\relax Cham,
  Switzerland: Springer Int. Pub., 2019, pp. 363--412.

\bibitem{Bioucas-Dias2012}
J.~Bioucas-Dias, A.~Plaza, N.~Dobigeon, M.~Parente, Q.~Du, P.~Gader, and
  J.~Chanussot, ``Hyperspectral unmixing overview: Geometrical, statistical,
  and sparse regression-based approaches,'' \emph{IEEE J. Sel. Topics Appl.
  Earth Observ. in Remote Sens.}, vol.~5, no.~2, pp. 354--379, 2012.

\bibitem{Bioucas2010}
J.~M. {Bioucas-Dias} and M.~A.~T. {Figueiredo}, ``Alternating direction
  algorithms for constrained sparse regression: {Application} to hyperspectral
  unmixing,'' in \emph{IEEE Works. Hypersp. Im. Sig. Proces.: Evolution Remote
  Sens.}, June 2010, pp. 1--4.

\bibitem{Rockafellar1970}
R.~Rockafellar, \emph{Convex Analysis}.\hskip 1em plus 0.5em minus 0.4em\relax
  New Jersey, USA: Princeton University Press, 1970.

\bibitem{Donoho1995}
D.~Donoho and I.~Johnstone, ``Adapting to unknown smoothness via wavelet
  shrinkage,'' \emph{J. Am. Statist. Assoc.}, vol.~90, no. 432, pp. 1200--1224,
  1995.

\bibitem{Combettes2011}
P.~L. Combettes and J.-C. Pesquet, ``Proximal splitting methods in signal
  processing,'' in \emph{Fixed-Point Algorithms for Inverse Problems in Science
  and Engineering}, H.~H. Bauschke, R.~S. Burachik, P.~L. Combettes, V.~Elser,
  D.~R. Luke, and H.~Wolkowicz, Eds.\hskip 1em plus 0.5em minus 0.4em\relax New
  York, NY: Springer, 2011, pp. 185--212.

\bibitem{Hintermuller2003}
M.~Hinterm\"{u}ller, K.~Ito, and K.~Kunisch, ``The primal--dual active set
  strategy as a semismooth {N}ewton method,'' \emph{SIAM J. Optim.}, vol.~13,
  no.~3, pp. 865--888, 2003.

\bibitem{Wen2012}
Z.~Wen, W.~Yin, H.~Zhang, and D.~Goldfarb, ``On the convergence of an
  active-set method for $\ell_1$ minimization,'' \emph{Optim. Meth. Soft.},
  vol.~27, no.~6, pp. 1127--1146, 2012.

\bibitem{Milzarek2014}
A.~Milzarek and M.~Ulbrich, ``A semismooth {Newton} method with
  multidimensional filter globalization for $\ell_1$-optimization,'' \emph{SIAM
  J. Optim.}, vol.~24, no.~1, pp. 298--333, 2014.

\bibitem{Solntsev2015}
S.~Solntsev, J.~Nocedal, and R.~H. Byrd, ``An algorithm for quadratic
  $\ell_1$-regularized optimization with a flexible active-set strategy,''
  \emph{Optim. Meth. Soft.}, vol.~30, no.~6, pp. 1213--1237, 2015.

\bibitem{DeSantis2016}
M.~De~Santis, S.~Lucidi, and F.~Rinaldi, ``A fast active set block coordinate
  descent algorithm for $\ell_1$-regularized least squares,'' \emph{SIAM J.
  Optim.}, vol.~26, no.~1, pp. 781--809, 2016.

\bibitem{Liang2017}
J.~Liang, J.~Fadili, and G.~Peyré, ``Activity identification and local linear
  convergence of forward--backward-type methods,'' \emph{SIAM J. Optimi.},
  vol.~27, no.~1, pp. 408--437, 2017.

\bibitem{Li2018}
X.~Li, D.~Sun, and K.~Toh, ``A highly efficient semismooth {Newton} augmented
  {Lagrangian} method for solving lasso problems,'' \emph{SIAM J. Optim.},
  vol.~28, no.~1, pp. 433--458, 2018.

\bibitem{Milzarek2019}
A.~Milzarek, X.~Xiao, S.~Cen, Z.~Wen, and M.~Ulbrich, ``A stochastic semismooth
  {Newton} method for nonsmooth nonconvex optimization,'' \emph{SIAM J.
  Optim.}, vol.~29, no.~4, pp. 2916--2948, 2019.

\bibitem{Facchinei2003a}
F.~Facchinei and J.-S. Pang, \emph{Finite-Dimensional Variational Inequalities
  and Complementarity Problems, Vol. I}.\hskip 1em plus 0.5em minus 0.4em\relax
  Springer-Verlag, 2003.

\bibitem{Qi1993b}
L.~Qi, ``Convergence analysis of some algorithms for solving nonsmooth
  equations,'' \emph{Math. Oper. Res.}, vol.~18, no.~1, pp. 227--244, 1993.

\bibitem{Ulbrich2011}
M.~Ulbrich, \emph{Semismooth {Newton} Methods for Variational Inequalities and
  Constrained Optimization Problems in Function Spaces}.\hskip 1em plus 0.5em
  minus 0.4em\relax Philadelphia, PA: MOS-SIAM Ser. Optim., 2011.

\bibitem{Qi1999}
L.~Qi and D.~Sun, \emph{A Survey of Some Nonsmooth Equations and Smoothing
  Newton Methods}.\hskip 1em plus 0.5em minus 0.4em\relax Boston, MA: Springer,
  1999, pp. 121--146.

\bibitem{Hintermuller2010}
\BIBentryALTinterwordspacing
M.~Hinterm\"{u}ller, ``Semismooth {N}ewton methods and applications,''
  Department of Mathematics, Humboldt-University of Berlin, Tech. Rep.
  November, 2010. [Online]. Available:
  \url{http://www.math.uni-hamburg.de/home/hinze/Psfiles/Hintermueller_OWNotes.pdf}
\BIBentrySTDinterwordspacing

\bibitem{Ulbrich2009}
M.~Ulbrich, ``Optimization methods in {Banach} spaces,'' in \emph{Optimization
  with PDE Constraints}.\hskip 1em plus 0.5em minus 0.4em\relax Dordrecht,
  Netherlands: Springer, 2009, pp. 97--156.

\bibitem{Vu2011}
B.~C. V{\~{u}}, ``A splitting algorithm for dual monotone inclusions involving
  cocoercive operators,'' \emph{Adv. Comput. Math.}, vol.~38, no.~3, pp.
  667--681, 2011.

\end{thebibliography}

\acrodef{ADMM}{alternating-direction method of multipliers}
\acrodef{BCCB}{block-circulant-circulant-block}
\acrodef{BSNR}{blurred-signal-to-noise ratio}
\acrodef{i.i.d.}{independent and identically distributed}
\acrodef{KKT}{Karush-Kuhn-Tucker}
\acrodef{PD}{positive-definite}
\acrodef{PSD}{positive semidefinite}	
\acrodef{VMPD}{variable-metric primal--dual method}
\acrodef{FFT}{fast Fourier Transform}
\acrodef{SNR}{signal-to-noise ratio}
\acrodef{RMSE}{root-mean-squared error}
\acrodef{CM}{method by Condat}		
\acrodef{LASSO}{least absolute shrinkage and selection operator}
\acrodef{KM}{Krasnosel'ski\u{\i}--Mann}
\acrodef{FB}{forward--backward}

\vfill

\end{document}